\documentclass{gtmon_a}
\pdfoutput=1

\usepackage[all]{xy}
\usepackage{pinlabel}

\proceedingstitle{Groups, homotopy and configuration spaces (Tokyo
  2005)}
\conferencestart{5 July 2005}
\conferenceend{11 July 2005}
\conferencename{Groups, homotopy and configuration spaces, 
                in honour of Fred Cohen's 60th birthday}
\conferencelocation{University of Tokyo, Japan}

\editor{Norio Iwase}
\givenname{Norio}
\surname{Iwase}

\editor{Toshitake Kohno}
\givenname{Toshitake}
\surname{Kohno}

\editor{Ran Levi}
\givenname{Ran}
\surname{Levi}

\editor{Dai Tamaki}
\givenname{Dai}
\surname{Tamaki}

\editor{Jie Wu}
\givenname{Jie}
\surname{Wu}

\title{A family of embedding spaces}

\author{Ryan Budney}
\givenname{Ryan}
\surname{Budney}
\address{Mathematics and Statistics\\
University of Victoria\\\newline
PO Box 3045 STN CSC\\
Victoria\\
British Columbia\\
V8W 3P4\\
Canada}
\email{rybu@uvic.ca}
\urladdr{}

\volumenumber{13}
\issuenumber{}
\publicationyear{2008}
\papernumber{3}
\startpage{41}
\endpage{83}

\doi{}
\MR{}
\Zbl{}

\arxivreference{math.AT/0605069}


\received{7 May 2006}
\revised{24 June 2007}
\accepted{2 July 2007}
\published{22 February 2008}
\publishedonline{22 February 2008}
\proposed{}
\seconded{}
\corresponding{}
\version{}

\makeatletter
\def\cnewtheorem#1[#2]#3{\newtheorem{#1}{#3}[section]
\expandafter\let\csname c@#1\endcsname\c@thm}

\let\xysavmatrix\xymatrix
\def\xymatrix{\disablesubscriptcorrection\xysavmatrix}
\AtBeginDocument{\let\bar\wbar\let\tilde\wtilde\let\hat\what}
\let\upepsilon\varepsilon

\makeop{cemb}
\makeop{supp}
\makeop{img}
\makeop{Id}

\newtheorem{thm}{Theorem}[section]
\cnewtheorem{lem}[thm]{Lemma}
\cnewtheorem{cor}[thm]{Corollary}
\cnewtheorem{prop}[thm]{Proposition}

\cnewtheorem{clm}[thm]{Claim}

\theoremstyle{definition}
\cnewtheorem{defn}[thm]{Definition}
\cnewtheorem{eg}[thm]{Example}
\cnewtheorem{question}[thm]{Question}

\makeatother  
\makeautorefname{defn}{Definition}
\makeautorefname{eg}{Example}

\newcommand{\I}{{\bf{I}}}

\newcommand{\Real}{{\Bbb R}}
\newcommand{\Rational}{{\Bbb Q }}
\newcommand{\Rat}{{\Bbb Q }}

\newcommand{\Zed}{{\Bbb Z}}

\newcommand{\Diff}{{\mathrm{Diff}}}

\newcommand{\SO}{{\mathrm{SO}}}

\newcommand{\Prime}{{\mathcal{P}}}

\newcommand{\Cu}{{\mathcal{C}}}

\newcommand{\CAut}{{\mathrm{CAut}}}

\newcommand{\K}{{\mathcal{K}}}
\newcommand{\PE}{{\mathcal{P}}}
\newcommand{\EK}[1]{{\mathrm{EC}({#1})}}
\newcommand{\ED}[1]{{\mathrm{ED}({#1})}}
\newcommand{\PK}[1]{{\mathrm{PEC}({#1})}}
\newcommand{\EC}[2]{{{\mathrm{Emb}}(\Real^{#1} \times {#2},\Real^{#1} \times {#2})}}
\newcommand{\Emb}{{\mathrm{Emb}}}

\newcommand{\splice}{{\bowtie}}

\newcommand{\gr}{{\mathrm{gr}}}

\begin{document}

\begin{htmlabstract}
Let Emb(S<sup>j</sup>,S<sup>n</sup>) denote the space of
C<sup>&infin;</sup>&ndash;smooth embeddings of the j&ndash;sphere in
the n&ndash;sphere.  This paper considers homotopy-theoretic
properties of the family of spaces Emb(S<sup>j</sup>,S<sup>n</sup>)
for n &ge; j> 0.  There is a homotopy-equivalence of
Emb(S<sup>j</sup>,S<sup>n</sup>) with SO<sub>n+1</sub>
&times;<sub>SO<sub>n-j</sub></sub> K<sub>n,j</sub> where
K<sub>n,j</sub> is the space of embeddings of <b>R</b><sup>j</sup> in
<b>R</b><sup>n</sup> which are standard outside of a ball.  The main
results of this paper are that K<sub>n,j</sub> is
(2n-3j-4)&ndash;connected, the computation of &pi;<sub>2n-3j-3</sub>
K<sub>n,j</sub> together with a geometric interpretation of the
generators.  A graphing construction &Omega; K<sub>n-1,j-1</sub>
&rarr; K<sub>n,j</sub> is shown to induce an epimorphism on homotopy
groups up to dimension 2n-2j-5.  This gives a new proof of
Haefliger's theorem that &pi;<sub>0</sub>
Emb(S<sup>j</sup>,S<sup>n</sup>) is a group for n-j>2.  The proof
given is analogous to the proof that the braid group has inverses.
Relationship between the graphing construction and actions of operads
of cubes on embedding spaces are developed.  The paper ends with a
brief survey of what is known about the spaces K<sub>n,j</sub>,
focusing on issues related to iterated loop-space structures.
\end{htmlabstract}

\begin{abstract}
Let $\mathrm{Emb}(S^j,S^n)$ denote the space of $C^{\infty}$--smooth
embeddings of the $j$--sphere in the $n$--sphere.  This paper
considers homotopy-theoretic properties of the family of
spaces $\mathrm{Emb}(S^j,S^n)$ for $n \geq j> 0$.  There is
a homotopy-equivalence $\mathrm{Emb}(S^j,S^n) \simeq SO_{n+1}
\times_{SO_{n-j}} \mathcal{K}_{n,j}$ where $\mathcal{K}_{n,j}$
is the space of embeddings of $\mathbb{R}^j$ in $\mathbb{R}^n$ which are
standard outside of a ball.  The main results of this paper are
that $\mathcal{K}_{n,j}$ is $(2n{-}3j{-}4)$--connected, the computation of
$\pi_{2n-3j-3} \mathcal{K}_{n,j}$ together with a geometric interpretation
of the generators.  A graphing construction $\Omega \mathcal{K}_{n-1,j-1}
\to \mathcal{K}_{n,j}$ is shown to induce an epimorphism on homotopy
groups up to dimension $2n{-}2j{-}5$.  This gives a new proof of Haefliger's
theorem that $\pi_0 \mathrm{Emb}(S^j,S^n)$ is a group for $n{-}j>2$.
The proof given is analogous to the proof that the braid group has
inverses.  Relationship between the graphing construction and actions of
operads of cubes on embedding spaces are developed.  The paper ends with
a brief survey of what is known about the spaces $\mathcal{K}_{n,j}$,
focusing on issues related to iterated loop-space structures.
\end{abstract}

\begin{asciiabstract}
Let Emb(S^j,S^n) denote the space of C^infty-smooth embeddings of the
j-sphere in the n-sphere.  This paper considers homotopy-theoretic
properties of the family of spaces Emb(S^j,S^n) for n >= j > 0.  There
is a homotopy-equivalence of Emb(S^j,S^n) with SO_{n+1}
times_{SO_{n-j}} K_{n,j} where K_{n,j} is the space of embeddings of
R^j in R^n which are standard outside of a ball.  The main results of
this paper are that K_{n,j} is (2n-3j-4)-connected, the computation
of pi_{2n-3j-3} (K_{n,j}) together with a geometric interpretation of
the generators.  A graphing construction Omega K_{n-1,j-1} --> K_{n,j}
is shown to induce an epimorphism on homotopy groups up to dimension
2n-2j-5.  This gives a new proof of Haefliger's theorem that pi_0
(Emb(S^j,S^n)) is a group for n-j>2.  The proof given is analogous to
the proof that the braid group has inverses.  Relationship between the
graphing construction and actions of operads of cubes on embedding
spaces are developed.  The paper ends with a brief survey of what is
known about the spaces K_{n,j}, focusing on issues related to iterated
loop-space structures.
\end{asciiabstract}

\maketitle

\section{Introduction}\label{INTRODUCTION}

Haefliger proved that the isotopy classes of smooth embeddings
of $S^j$ in $S^n$ form a group provided $n-j>2$, with the connect-sum as
multiplication.   This paper starts with a new proof of Haefliger's
 result, showing not only that $\pi_0 \Emb(S^j,S^n)$ is a group, but the 
{\it reason} it is a group is that every element is {\it spun} (see
\fullref{lgh} for the definition of the graphing/spinning map, $\gr_1$). 
The inverse of a spun knot is its mirror-reflection, as in braid
groups.  The key strategy revolves around a pseudo-isotopy
fibre-sequence $\K_{n+1,j+1} \to \PE_{n,j} \to \K_{n,j}$. The fact that
the pseudo-isotopy embedding space $\PE_{n,j}$ is connected implies the result.
In his dissertation, Tom Goodwillie \cite{GoodD} gave a very detailed study 
of (general) pseudo-isotopy embedding spaces. His results include that 
$\PE_{n,j}$ is at least $(2n{-}2j{-}5)$--connected. This allows for the computation of
the first non-trivial homotopy groups of $\K_{n,j}$ and $\Emb(S^j,S^n)$ provided
$2n-3j-3 \geq 0$. The 2--fold spinning construction 
$\pi_2 \K_{4,1} \to \pi_0 \K_{6,3} = \pi_0 \Emb(S^3,S^6) \simeq \Zed$
is shown to be an isomorphism, answering a question posed by the author
in \cite{Cube}.  This also allows for a new construction of explicit
generators of $\pi_{2n-3j-3} \K_{n,j}$ for all $n,j$ such that $2n-3j-3
\geq 0$.


\begin{defn}\label{BIGdef}
\quad
\begin{itemize}
\item $D^n := \{ x \in \Real^n : |x|\leq 1 \}$ is the unit $n$--disc,
with $S^{n-1} = \partial D^n$ the $(n{-}1)$--sphere. 
\item $\I = [-1,1] = D^1$ is the standard interval.
\item Given a topological space (resp. smooth manifold) $X$ with base-point, 
denote the space of continuous (resp. smooth) functions
$f\co \Real \to X$ such that $f(\Real \setminus \I)=*$ by $\Omega X$. 
\item $\Emb(D^j,D^n)$ denotes the space of embeddings $f \co D^j \to D^n$
which are `neat' in the sense that $f(D^j) \cap S^{n-1} = f(S^{j-1})$ and 
$f$ intersects $S^{n-1}$ transversely.
\item The space of smooth embeddings of a $j$--sphere in 
an $n$--sphere is denoted $\Emb(S^j,S^n)$. 
\item $\K_{n,j}$ denotes the space of `long' embeddings of $\Real^j$ in 
$\Real^n$. This is the space of all smooth embeddings $f \co \Real^j \to \Real^n$ 
such that
$$f(t_1,t_2,\ldots,t_j)=(t_1,t_2,\ldots,t_j,0,\ldots,0)$$
provided $(t_1,\ldots,t_j) \not\in \I^j$ and 
$f(\Real^j) \cap \partial \I^n = \partial \I^j \times \{0\}^{n-j}$.
If $f \in \K_{n,j}$, let $\K_{n,j}(f)$ denote the path-component of
$\K_{n,j}$ containing $f$.  We will show $\K_{n,j}$ has the homotopy
type of the subspace of $\Emb(D^j,D^n)$ such that every embedding restricts
to a fixed linear embedding on the boundary. 
\item  Let $\PE_{n,j}$ denote the space of embeddings $f \co \Real^j \to \Real^n$ such that:
\begin{itemize} 
\item $f(t_1,t_2,\ldots,t_j) = (t_1,t_2,\ldots,t_j,0,\ldots,0)$ 
for $(t_1,\ldots,t_j) \not\in [-1,\infty)\times \I^{j-1}$ 
\item there is a $g \in \K_{n-1,j-1}$ such that for all
$(t_1,\ldots,t_j) \in [1,\infty)\times \Real^{j-1}$,  
$f(t_1,t_2,\ldots,t_j)=(t_1,g(t_2,\ldots,t_j))$.
\item $f(\Real^j) \cap \partial \I^n = f(\partial \I^j) \times \{0\}^{n-j}$.
\end{itemize}
In the literature, $\PE_{n,j}$ is sometimes 
given the notation $PE(D^{j-1},D^{n-1})$, $C(D^{j-1},$ $D^{n-1})$ or 
$\cemb(D^{j-1},D^{n-1})$, and is either called a pseudoisotopy embedding 
space, or concordance embedding space respectively.
Here it will be called the pseudoisotopy embedding space. We will show
that $\PE_{n,j}$ has the homotopy-type of the subspace of
$\Emb(D^j,D^n)$ which restricts to a standard linear embedding
on a fixed hemisphere in the boundary of $D^j$. 
\item $\EK{j,M}$ is defined to be the space of embeddings
$f \co  \Real^j \times M \to \Real^j \times M$ such that
$\supp(f) \subset \I^j \times M$,
where, $\supp(f) = \{ x \in \Real^j \times M : f(x)\neq x \}$. 
`EC' stands for `cubically-supported embeddings'. We are mostly
interested in the case where $M$ is a disc $M=D^k$.
These embeddings are not required to send boundary to boundary.
See \fullref{fig1}.
\begin{figure}\label{fig1}
\centering
\labellist\small
\pinlabel {$f\in\EK{1,D^2}$} [t] at 400 378
\pinlabel {$1$} [t] at 670 20
\pinlabel {$-1$} [t] at 105 20
\endlabellist
$$\includegraphics[width=10cm]{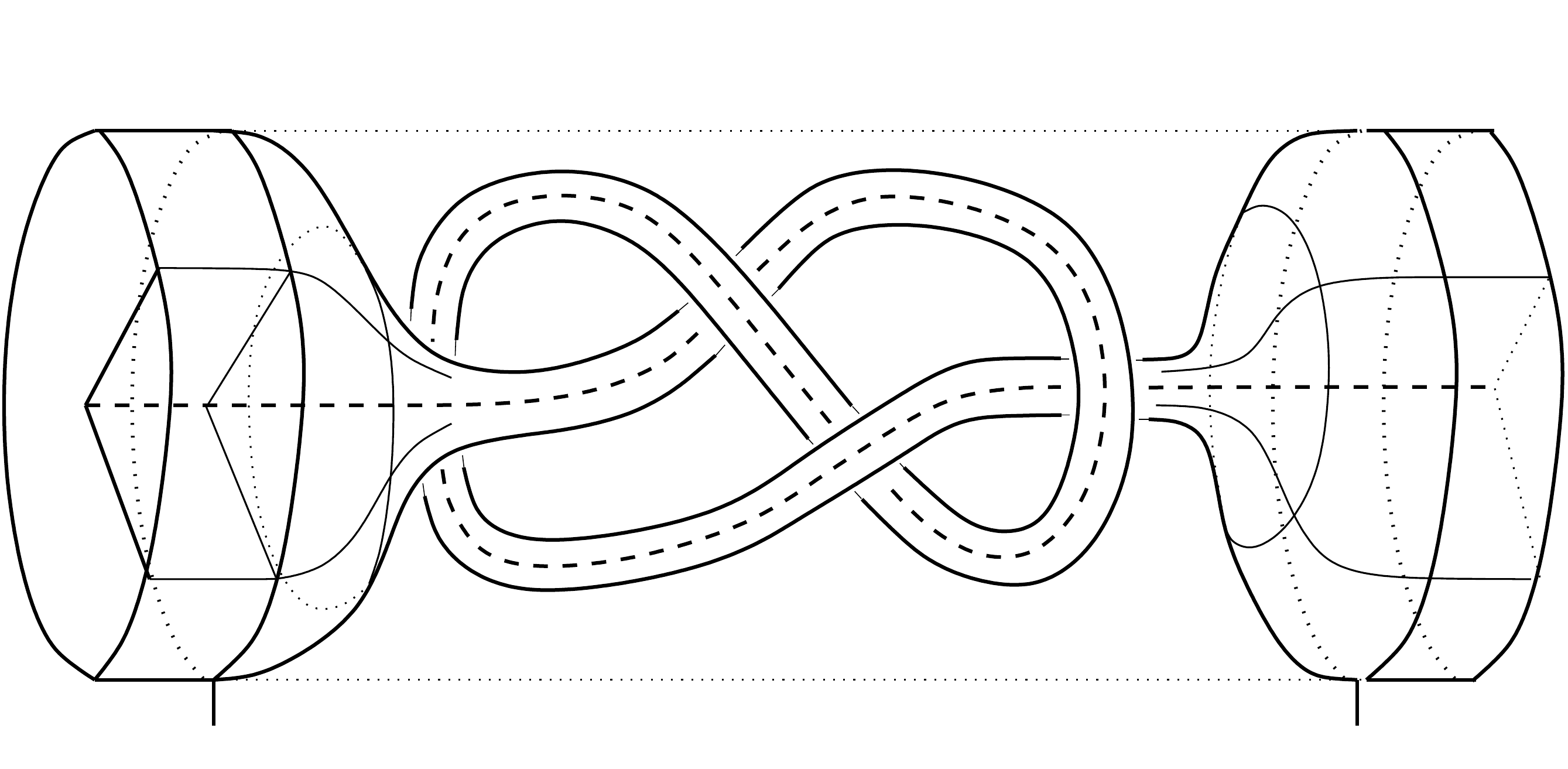}$$
\caption{}
\end{figure}
\item $\PK{j,M}$ is the space of embeddings $f \co  \Real^j \times M \to \Real^j \times M$
such that $\supp(f) \subset [-1,\infty) \times \I^{j-1} \times M$
and there exists some $g \in \EK{j-1,M}$ such that 
$f(t_1,t_2,\ldots,t_j,m) = (t_1,g(t_2,\ldots,t_j,m))$ for all $(t_1,t_2,\ldots,t_j,m) \in [1,\infty) \times \Real^{j-1} \times M$.
The letters `PEC' stand for `cubically-supported embedding pseudo-isotopy space.'
$\PK{j,D^k}$ is the framed analogue of $\PE_{j+k,j}$.
\item  A diagram of two maps $A \to C \to D$ is a homotopy
fibre sequence if there exists a commutative diagram
$$\xymatrix{A \ar[r] \ar[d] & C \ar[r] \ar[d] & D \ar[d] \\
            F \ar[r] & E \ar[r] & B }$$
such that $F \to E \to B$ is a fibration and the vertical maps
are homotopy-equivalences. 
\item $\Diff(D^n)$ denotes the space of smooth diffeomorphisms of $D^n$
which restrict to the identity on the boundary.  $\Diff(S^n)$ is the group
of diffeomorphisms of $S^n$.
\end{itemize}
\end{defn}

All embedding spaces are endowed with the weak $C^\infty$--topology 
(see Hirsch \cite{Hirsch}), sometimes also called the Whitney topology. Many classical 
results on the homotopy properties of embedding spaces that will be 
repeatedly used in this paper appear in Cerf's paper \cite{Cerf}, such
as the fibration properties of restriction maps, and the homotopy-classification
of spaces of tubular neighbourhoods. 

In the definition of $\K_{n,j}$ replacing the cubes $\I^n$ and $\I^j$ with
discs $D^n$ and $D^j$ gives a homotopy-equivalent space. Similarly for the
definition of $\Diff(D^n)$ and $\EK{j,M}$. The proof is a typical argument
when one deals with these spaces, see for example \cite[Corollary~6]{Cube}.
 
\fullref{fibrations} briefly covers the most elementary
relationships between the spaces defined above: 
$\K_{n,j}$, $\Emb(S^j,S^n)$, $\Emb(S^j,\Real^n)$, $\Emb(D^j,D^n)$, 
$\PE_{n,j}$, $\EK{j,D^{n-j}}$ and $\PK{j,D^{n-j}}$. This
section also includes a generalisation of an observation of
Goodwillie and Sinha \cite{Dev2} concerning the Smale--Hirsch
map $\K_{n,j} \to \Omega^j V_{n,j}$.  The Goodwillie--Sinha
result is that this map is null-homotopic for $j=1$. The
generalisation that appears here is that the map factors 
as a composite
$\K_{n,j} \to \Omega^j V_{n-1,j-1} \to \Omega^j V_{n,j}$ where 
the map $\Omega^j V_{n-1,j-1} \to \Omega^j V_{n,j}$ is the $j$--fold
loop of the fibre inclusion in the Stiefel fibration
$V_{n-1,j-1} \to V_{n,j} \to S^{n-1}$.

\fullref{path} is the heart of the paper.  A proof of
Haefliger's theorem, that for $n-j>2$ $\pi_0 \Emb(S^j,S^n)$ is a group is given.  
The proof permutes some of the main concepts of Haefliger's original argument. 
It has two essential steps:
1) The construction of a homotopy-equivalence
$\Emb(S^j,S^n) \simeq \SO_{n+1} \times_{\SO_{n-j}} \K_{n,j}$ together with
fibrations $\PE_{n,j} \to \Emb(D^j,D^n) \to V_{n,j}$ and
$\K_{n,j} \to \PE_{n,j} \to \K_{n-1,j-1}$
reduces the problem to 2) proving that $\Emb(D^j,D^n)$ is connected.
Thus, the argument boils down to showing the monoid $\pi_0 \K_{n,j}$ 
is a group because it is the image of the group $\pi_1 \K_{n-1,j-1}$. Further, it is 
shown that the `boundary map' $\gr_1 \co  \Omega \K_{n-1,j-1} \to \K_{n,j}$ has a geometric interpretation as a variant of Litherland `deform spinning.' In this case it is given by
the formula
$$(\gr_1 f)(t_0,t_1,\ldots,t_{j-1})
= \left(t_0,f(t_0)\left(t_1,\ldots,t_{j-1}\right)\right).$$
In \fullref{Connectivityknj}, Goodwillie's dissertation is used
to prove that $\gr_1 \co  \Omega \K_{n-1,j-1} \to \K_{n,j}$ induces an epimorphism of the 
on homotopy groups $\pi_i$ for $i \leq 2n-2j-5$. By comparing with the work of
Turchin and Sinha this allows the computation of 
$\pi_{2n-3j-3} \K_{n,j}$.  An enumerative-geometry argument is used to
construct a cohomology class $\nu_2 \in H^{2n-6}(\K_{n,1};\Zed)$, which
is used to find an explicit generator of $\pi_{2n-6} \K_{n,1} \simeq \Zed$.
The generator can be thought of as the resolutions of a long immersion of 
$\Real$ in $\Real^n$ having two regular double points, corresponding
to the $\bigotimes$ chord diagram.  The generators of the groups $\pi_0 \K_{n,j}$ 
for $2n-3j-3=0$ are constructed as iterated graphs of the generator of 
$\pi_{j-1} \K_{n-j+1,1}$.

\fullref{actions} investigates the extent to which
the fibration $\K_{n,j} \to \PE_{n,j} \to \K_{n-1,j-1}$, and its
framed analogue are equivariant with respect to natural actions of 
operads of cubes.  $\PK{j,M}$ is shown to have an
action of the operad of $j$--cubes, the map $\PK{j,M} \to \EK{j-1,M}$
is shown to be equivariant with respect to the $j$--cubes action
defined by the author in \cite{Cube}.  The graphing construction 
$\Omega \EK{j-1,M} \to \EK{j,M}$ is shown to be equivariant with respect 
to the $(j+1)$--cubes action.

\fullref{survey} covers, in a rather terse survey manner, many of the
basic properties the spaces $\K_{n,j}$ which have not already been mentioned.
A curiosity is put forward: two seemingly distinct null homotopies of the 
inclusion $\K_{n,1} \to \K_{n+1,1}$ are described, giving a mysterious map 
$\Sigma \K_{n,1} \to \K_{n+1,1}$. This leads to a question about the existence
of a `Freudenthal suspension' $\Sigma^2 \K_{n,1} \to \K_{n+1,1}$. Basic properties 
of other natural maps such as $\K_{n,j} \to \Omega \K_{n,j-1}$ and the
Smale--Hirsch map 
$SH \co  \K_{n,j} \to \Omega^j V_{n,j}$ are described.  

Part of this manuscript was produced while visiting the University
of Rome `La Sapienza', Louvain-la-neuve, the American Institute of
Mathematics, the University of Tokyo and IH\'ES. I would especially like to
thank the Max Planck Institute for Mathematics, in Bonn, for giving me
the freedom to pursue this line of enquiry.
I would like to thank my hosts for their hospitality: Riccardo Longoni, 
Paolo Salvatore,  Corrado De Concini, Magnus Jacobsson, Pascal Lambrechts, 
Victor Turchin, and Toshitake Kohno. Victor Turchin's comments on the first 
draft of this manuscript were particularly helpful. I would like to thank
several mathematicians whose comments, knowingly or not, have helped me in putting
this paper together: Greg Arone, John Rognes, Tom Goodwillie, Larry Siebenmann,
Dev Sinha, Arkadiy Skopenkov, Lee Rudolph, Matthias Kreck, Paolo Salvatore,
Jianguo Cao and Danny Ruberman.

\section{Basic relations between embedding spaces}\label{fibrations}

This section describes some basic relationships between the spaces: 
$\K_{n,j}$, $\EK{j,M}$, $\Emb(S^j,S^n)$, $\Emb(S^j,\Real^n)$,
$\Emb(D^j,D^n)$, $\PE_{n,j}$ and $\PK{j,M}$.  The
essential spirit of the results is that most homotopy 
questions about these spaces reduce to studying $\K_{n,j}$ and $\PE_{n,j}$.

Given a neat embedding  $f \co  D^j \to D^n$, the restriction to the 
boundary is an embedding $f_{|\partial D^j} \co  S^{j-1} \to S^{n-1}$. On a 
global level, restriction defines a function
$$ \Emb(D^j,D^n) \to \Emb(S^{j-1},S^{n-1}) $$ 
which is a fibration (see Cerf \cite{Cerf} and Palais \cite{Pal}).  In this paper `fibration'
means Serre fibration. The above map is known to be more than a fibration,
it is a locally trivial fibre-bundle \cite{Pal}. Fibrations
need not be onto.  In this example, the fibration is onto the isotopy classes
of `slice' knots (and not all knots are slice, see Kawauchi \cite{Kaw} for 
examples).  Thus, the homotopy-type of the fibre can change as one changes 
base-space components, and fibres are allowed to be empty.

Consider $\Emb(S^{j-1},S^{n-1})$ to be a based space, with base-point the 
standard inclusion $S^{j-1} \subset S^{n-1}$. The fibre of 
$\Emb(D^j,D^n) \to \Emb(S^{j-1},S^{n-1})$
over the base-point has the homotopy-type of $\K_{n,j}$. 
There is a similar fibration $\K_{n,j} \to \PE_{n,j} \to \K_{n-1,j-1}$
defined by restriction to the `free face.' The next theorem shows that
this fibration induces the fibration $\Emb(D^j,D^n) \to \Emb(S^{j-1},S^{n-1})$.

\begin{thm}\label{Splitprop}
For $n-j>0$ there are homotopy-equivalences
\begin{align*}
\Emb(D^j,D^n) &\simeq \SO_n \times_{\SO_{n-j}} \PE_{n,j} \\
\Emb(S^{j-1},S^{n-1}) &\simeq \SO_n \times_{\SO_{n-j}} \K_{n-1,j-1}.
\end{align*}
Moreover, the homotopy fibre sequence 
$\K_{n,j} \to \Emb(D^j,D^n) \to \Emb(S^{j-1},S^{n-1})$
fits into a commutative diagram of $6$ homotopy fibre sequences:

$$\xymatrix{
\K_{n,j} \ar[r] \ar[d] & \PE_{n,j} \ar[r] \ar[d]     & \K_{n-1,j-1}  \ar[d] \\
\K_{n,j} \ar[r] \ar[d] & \Emb(D^j,D^n) \ar[r] \ar[d] & \Emb(S^{j-1},S^{n-1}) \ar[d] \\
\text{*} \ar[r]               & V_{n,j} \ar[r]              & V_{n,j} 
}$$
\end{thm}
\begin{proof}
In Budney and Cohen \cite{BudCoh} a homotopy-equivalence 
$$\SO_n \times_{\SO_{n-j}} \K_{n-1,j-1} \to \Emb(S^{j-1},S^{n-1})$$
was constructed. The basic
idea is to consider $S^{n-1}$ to be the one-point compactification of
$\Real^{n-1}$, this gives an inclusion $\K_{n-1,j-1} \to \Emb(S^{j-1},S^{n-1})$.
The action of $\SO_n$ on $S^{n-1}$ gives an extension
$$\SO_n \times_{\SO_{n-j}} \K_{n-1,j-1} \to \Emb(S^{j-1},S^{n-1}).$$ 
$\SO_n \times_{\SO_{n-j}} \K_{n-1,j-1}$ fibres over $V_{n,j} = \SO_n / \SO_{n-j}$
by projection onto the first coordinate.  $\Emb(S^{j-1},S^{n-1})$ fibres over
a space homotopy-equivalent to $V_{n,j}$ by restriction to a fixed hemi-sphere
$B \subset S^{j-1}$, $\Emb(S^{j-1},S^{n-1}) \to \Emb(B,S^{n-1}) \simeq V_{n,j}$ 
\cite{Cerf}. This makes $\SO_n \times_{\SO_{n-j}} \K_{n-1,j-1} \to \Emb(S^{j-1},S^{n-1})$
a map of fibrations.

The same idea can be applied to $\Emb(D^j,D^n)$. Let $B \subset \partial D^j = S^{j-1}$ 
be as above.  Let $\Emb(D^j \text{ rel } B,D^n)$ denote the
subspace of $\Emb(D^j,D^n)$ which is fixed point-wise on $B$.
There is a fibre bundle 
$\Emb(D^j \text{ rel } B,D^n) \to \Emb(D^j,D^n) \to \Emb( B,S^{n-1})$ 
given by restriction to $B$.
The base-space has the homotopy-type of $V_{n,j} \simeq \SO_n/\SO_{n-j}$
and as in the previous paragraph, there is a map of fibrations
$$\SO_n \times_{\SO_{n-j}} \Emb(D^j \text{ rel } B,D^n) \to \Emb(D^j,D^n).$$
That $\Emb(D^j \text{ rel } B,D^n)$ has the same homotopy-type as 
$\PE_{n,j}$ is a fairly standard argument, see for example the second half of 
\cite[Corollary~6]{Cube}. 
\end{proof}

When $n=j$, the above argument proves that $\Emb(D^n,D^n)$
has the homotopy-type of  $O_n \times \PE_{n,n}$. Similarly,
$\Emb(S^{n-1},S^{n-1}) = \Diff(S^{n-1})$ has the homotopy-type of $O_n
\times \K_{n-1,n-1}$. This case appears in Hatcher \cite{Hatcher2}.

There is a similar relationship between $\Emb(S^j,\Real^n)$ and
$\K_{n,j}$. For this proposition, identify 
$\overline{\Real^n}$ (the one-point compactification of $\Real^n$) with $S^n$
via stereographic projection.  This makes $\SO_n$ the stabiliser of $\infty$
under the $\SO_{n+1}$ action on $S^n$.  Denote the projection map
$\SO_{n+1} \to S^n$ by $\pi$.  Given $f \in \K_{n,j}$ let $\wbar{f} \in \Emb(S^j,S^n)$ be the one-point compactification of $f$. Notice that the space
$$ \{(A,f) : A \in \SO_{n+1}, \pi(A) \in S^n \setminus \img(\wbar{f}), f \in \K_{n,j} \}$$
fibres over $C \rtimes \K_{n,j}$ with fibre $\SO_n$, for
$$C \rtimes \K_{n,j} = \{(p,f) : p \in S^n \setminus \img(\wbar{f}), f \in \K_{n,j} \}.$$
Denote $\{(A,f) : A \in \SO_{n+1}, \pi(A) \in S^n \setminus \img(\wbar{f}), f \in \K_{n,j} \}$
by $(C \rtimes \K_{n,j})^*(\pi)$. 
Consider $(C \rtimes \K_{n,j})^*(\pi)$ to be the pull-back of $\pi$ over $\Real^n$.  
Since $\pi$ is trivial over $\Real^n$, the pull-back must be as well.
$$ \SO_n \times (C \rtimes \K_{n,j}) \simeq (C \rtimes \K_{n,j})^*(\pi).$$
 Notice that $\SO_{n-j}$ acts on $(C \rtimes \K_{n,j})^*(\pi)$ from the left, 
by considering $\SO_{n-j} \subset \SO_{n+1}$
to be the group that leaves $S^j = \overline{\Real^j}$ in $S^n$ fixed point-wise.

\begin{prop}\label{long-euclid}
Provided $n-j>0$ there is a homotopy-equivalence 
$$
\SO_{n-j} \backslash (C \rtimes \K_{n,j})^*(\pi) \to \Emb(S^j,\Real^n)
$$
induced by the map $(A,f) \longmapsto A^{-1} \circ \wbar{f}$.
Moreover, there is a homotopy-equivalence
$$\SO_{n-j} \backslash (C \rtimes \K_{n,j})^*(\pi) \to \SO_n \times_{\SO_{n-j}} (C \rtimes \K_{n,j})$$
where the action of $\SO_{n-j}$ on $\SO_n$ is by left multiplication.
\end{prop}
\begin{proof}
Observe that $\Emb(S^j,\Real^n)$ fibres over $V_{n,j}$.  The fibre
can be identified with $\{ f \in \K_{n,j} : 0 \not\in f(\Real^j) \}$. $C \rtimes \K_{n,j}$
fibres over a ball with fibre $\{ f \in \K_{n,j} : 0 \not\in f(\Real^j) \}$, thus there is
a homotopy-fibre sequence
$$ C \rtimes \K_{n,j} \to \Emb(S^j,\Real^n) \to V_{n,j}$$
$ (C \rtimes \K_{n,j})^*(\pi)$ similarly fibres over $V_{n,j}$ giving a commutative
ladder of homotopy fibre sequences
$$\xymatrix{
C \rtimes \K_{n,j} \ar[r] & \Emb(S^j,\Real^n) \ar[r] & V_{n,j} \\
C \rtimes \K_{n,j} \ar[r] \ar[u] & (C \rtimes \K_{n,j})^*(\pi) \ar[r] \ar[u] & V_{n,j} \ar[u]
}
$$
Let $(A,f) \in (C \rtimes \K_{n,j})^*(\pi)$, then $A$ is a matrix whose first column
vector is $\pi(A)$, the remaining vectors are in the tangent space to $\Real^n$ at
$\pi(A)$.  Let $[A]_{\pi(A)}$ denote the representation of $A$ with respect to the
standard framing of $\Real^n$ at $\pi(A)$. Consider the map 
$(C \rtimes \K_{n,j})^*(\pi) \to SO_n \times (C \rtimes \K_{n,j})$
given by sending the pair $(A,f)$ to $\left([A]_{\pi(A)}, (\pi(A), f)\right)$. 
This map is equivariant with respect to the action of $\SO_{n-j}$
since if $B \in \SO_{n-j}$ then $B.(A,f) = (BA,Bf)$, which is sent to
$\left([BA]_{\pi(BA)}, (\pi(BA),Bf)\right) = \left( [BA]_{B\pi(A)}, B.(A,f)\right)$, 
but $[BA]_{B\pi(A)} = B[A]_{\pi(A)}$ by a change of variables argument, giving the result.
\end{proof}

A basic fact and conventions about homotopy-fibres is given for future
reference.

\begin{lem}\label{hom-fib}
Let $p \co  E \to B$ be a fibration. Let $e \in E$ and $b \in B$ be 
the base-points of $E$ and $B$ respectively, with $p(e)=b$.
Let $i \co  F \to E$ be the fibre inclusion. Let $R(F) = \{(a,h) : a \in F, 
h \co  [0,1] \to E, h(0)=i(a) \}$ then the map
$R(i) \co  R(F) \to E$ given by evaluation $h(1)$ is a fibration,
and $\pi_F \co  R(F) \to F$ given by projection onto $F$ is a
homotopy-equivalence. The fibre of the map $R(i) \co  R(F) \to E$
is the space $HF(i) = \{h : [0,1] \to E, h(0) \in F, h(1)=e\}$, and
the map $p_* \co  HF(i) \to \Omega B$ given by post-composition with $p$
is a weak homotopy-equivalence, giving a fibration:
$$\Omega E \to HF(i) \to F$$
and a homotopy-commutative diagram
$$\xymatrix{
                & \Omega B      & F \ar[r]^i & E \ar[r]^p & B \\
\Omega E \ar[r] & HF(i) \ar[u]^-{p_*} \ar[r] \ar[ur] & R(F) \ar[u]^{\simeq}_{\pi_F} \ar[ur]_{R(i)} &   
}$$
\end{lem}

The map $HF(i) \to F$ is sometimes called the `connecting map' or the
`boundary map' as it induces the same map as the connecting map in the
homotopy long exact sequence of the fibration $p$. 

The next two results are a modest generalisation of observations due to 
Goodwillie (unpublished), Sinha \cite{Dev2}, Turchin \cite{Tur3} and 
Salvatore \cite{Salvatore}, concerning the monodromy of the
fibration $\EK{j,D^{n-j}} \to \K_{n,j}$ and the Smale--Hirsch map
$\K_{n,j} \to \Omega^j V_{n,j}$. Note, $\Omega^j V_{n,j}$ has the
homotopy-type of the space of long immersions $\Real^j \to \Real^n$
provided $n-j>0$, by the Smale--Hirsch theorem. 

\begin{thm}\label{Trivprop}
The homotopy fibre sequence
$$\Omega^j \SO_{n-j} \to \EK{j,D^{n-j}} \to \K_{n,j}$$
is trivial for $j=1$, and also for $n-j \leq 2$.
There is a pull-back diagram
of homotopy fibre sequences:
$$\xymatrix{\Omega^j\SO_{n-j} \ar[d] \ar[r] & \Omega^j \SO_{n-j} \ar[d] \\
            \EK{j,D^{n-j}} \ar[d] \ar[r]   & P \Omega^{j-1} \SO_{n-j} \ar[d] \\
            \K_{n,j} \ar[r]^-{cl}          & \Omega^{j-1} \SO_{n-j} 
}$$
Where $\Omega^j \SO_{n-j} \to P\Omega^{j-1} \SO_{n-j} \to \Omega^{j-1} \SO_{n-j}$ is the 
path-loop fibration of the space $\Omega^{j-1} \SO_{n-j}$. 
The classifying map $cl\co \K_{n,j} \to \Omega^{j-1} \SO_{n-j}$ fits into
a commutative diagram
$$\xymatrix{ 
 & & \Omega^j V_{n,n-j} \ar[d] \\
\Omega^j \SO_j \ar[r] & \Omega^j V_{n,j} \ar[r] & \Omega^j G_{n,j} \equiv \Omega^j G_{n,n-j} \ar[d]^{\text{mono}} &  \\
 & \K_{n,j} \ar[u]^{SH} \ar[r]^-{cl} & \Omega^{j-1} \SO_{n-j} }$$
where `$SH$' is the Smale--Hirsch map,
$V_{n,j}$ is the Stiefel manifold of $j$ linearly independent vectors in 
$\Real^{n}$, $\SO_{j} \to V_{n,j} \to G_{n,j}$ is the canonical fibration
for the Grassmanian of oriented $j$--dimensional subspaces of $\Real^n$.
`mono' is the $j$--fold looping of the classifying map 
$G_{n,n-j} \to B\SO_{n-j}$ for the bundle $\SO_{n-j} \to V_{n,n-j} \to G_{n,n-j}$. Identify
$G_{n,j}$ with $G_{n,n-j}$ via the oriented orthogonal complement.

Framed and unframed pseudoisotopy embedding spaces are more
directly related, as the forgetful map
$\PK{j,D^{n-j}} \to \PE_{n,j}$ is a homotopy-equivalence.
\end{thm}

\begin{proof} 
The observation of the existence of the above pull-back 
diagram first appears in Turchin's work \cite{Tur3} for $j=1$. The idea
is to divide $\I^j$ into $\I \times \I^{j-1}$. Given a knot $f \in \K_{n,j}$,
let $\nu f$ be its normal bundle, and consider parallel transport (using
the connection inherited as a submanifold of Euclidean space $\Real^n$) from
$\nu f_{|\{-1\} \times \I^{j-1}}$ to 
$\nu f_{\{1\} \times \I^{j-1}}$, this is an element of
$\Omega^{j-1} \SO_{n-j}$.  The map $\EK{j,D^{n-j}} \to P\Omega^{j-1}\SO_{n-j}$
is defined similarly, only along the paths $\I \times \{x\} \subset \I \times \I^{j-1}$
$f \in \EK{j,D^{n-j}}$ one has a pre-defined framing of $\nu f_{|\Real^j\times\{0\}^{n-j}}$
which can be compared to the parallel transport framing, giving the bundle map.  

Observe that the way 
$\K_{n,j} \to \Omega^{j-1} \SO_{n-j}$ is defined, it factors as a composite
$\K_{n,j} \to \Omega^j G_{n,j} \equiv \Omega^j G_{n,n-j} \to \Omega^{j-1} \SO_{n-j}$.
$\K_{n,j} \to \Omega^j G_{n,j}$ is the `tangent space map.' $G_{n,j}$ is the 
Grassmanian of $j$--dimensional subspaces of $\Real^n$.
$\text{mono} \co  \Omega^j G_{n,n-j} \to \Omega^{j-1} \SO_{n-j}$ is the $j$--fold
looping of the classifying map of the bundle $\SO_{n-j} \to V_{n,n-j} \to G_{n,n-j}$.

For the fibration $\PK{j,D^{n-j}} \to \PE_{n,j}$ observe the fibre has the
homotopy-type of the path-space $P\Omega^{j-1} \SO_{n-j}$.
\end{proof}

The homotopy-class of the Smale--Hirsch map
$SH \co  \K_{n,j} \to \Omega^j V_{n,j}$ is not so well understood.
There are results concerning the induced map $SH \co  \pi_0 \K_{n,j} \to \pi_j V_{n,j}$ in two 
cases: Kervaire proved it to be trivial provided $2n-3j \geq 2$ \cite{Ker}.  
In the co-dimension 2 case $n-j=2$, Hughes and Melvin showed that
$SH \co  \pi_0 \K_{n,j} \to \pi_j V_{n,j}$ has non-trivial image if and only
if $j \equiv 3 \text{ mod } 4$ \cite{HughMelv}, moreover they gave a rather appealing
description of the immersions that can be realised as embeddings. 
Eckholm and Sz\"ucs \cite{Ek,EkSz} have recently given more geometric
interpretations of the obstruction to an immersion having a representative
that is an embedding.

\begin{thm}\label{gsproof}
The Smale--Hirsch map $SH\co  \K_{n,j} \to \Omega^j V_{n,j}$
fits into a homotopy-commutative diagram
$$\xymatrix{ 
\K_{n,j} \ar[dr] \ar[rr]^{SH} & & \Omega^j V_{n,j} \\
 & \Omega^j V_{n-1,j-1} \ar[ur]_-{\Omega^j(i)} &
}$$
where $i \co  V_{n-1,j-1} \to V_{n,j}$ is the fibre-inclusion
of the fibration $V_{n-1,j-1} \to V_{n,j} \to S^{n-1}$.
\begin{proof}
Consider the commutative diagram of spaces and maps:
$$\xymatrix{
\K_{n,j} \ar[r] \ar[d]^{SH} & \PE_{n,j} \ar[r] \ar[d]^{SH} & \K_{n-1,j-1} \ar[d]^{SH} \\
\Omega^j V_{n,j} \ar[r] & \Omega^{j-1} HF(i) \ar[r] & \Omega^{j-1} V_{n-1,j-1}
}$$
$HF(i)$ is the homotopy-fibre of $i$. 
By \fullref{hom-fib}, there is a homotopy-equivalence 
$HF(i) \simeq \Omega S^{n-1}$.

The Smale--Hirsch map $SH \co  \PE_{n,j} \to \Omega^j S^{n-1}$ is given by 
differentiation in the vertical `pseudo-isotopy' direction. 
The map $h\co  [0,3]\times \Real^j \times \PE_{n,j} \to S^{n-1}$ given by
$$
h(t,x_1,\ldots,x_j,f) =
\left\{
\begin{array}{ll}
n\bigl(\frac{\partial f}{\partial x_1}(x_1,\ldots,x_j)\bigr) & t=0 \\
n(f(x_1+t,x_2,\ldots,x_j)-f(x_1,\ldots,x_j))              & 0<t\leq 2 \\
p_{t-2}(n(f(x_1+2,x_2,\ldots,x_j)-f(x_1,\ldots,x_j)))     & 2 \leq t \leq 3
\end{array}
\right.
$$
is a null-homotopy of the Smale--Hirsch map, 
provided $p \co  [0,1] \times S^{n-1} \setminus \{-1\} \to S^{n-1} \setminus \{-1\}$
is a deformation-retraction of 
$S^{n-1} \setminus \{-1\}$ to $\{1\} \subset S^{n-1}$,
and $n \co  \Real^n \setminus \{0\} \to S^{n-1}$ is the function
$n(v)=\frac{v}{|v|}$. 
\end{proof}
\end{thm}

Theorems \ref{Trivprop} and \ref{gsproof} combine to give
a commutative diagram involving the maps $cl \co  \K_{n,j} \to \Omega^{j-1} \SO_{n-j}$
and $SH \co  \K_{n,j} \to \Omega^j V_{n,j}$.
$$\xymatrix{ 
 & & \Omega^j V_{n-1,n-j} \ar[d] \\
\Omega^j V_{n,j} & \Omega^j V_{n-1,j-1} \ar[l] \ar[r]^-{\Omega^j\perp} & 
\Omega^j G_{n-1,j-1} \equiv \Omega^j G_{n-1,n-j} \ar[d] \\
 & \K_{n,j} \ar[ul]^{SH} \ar[u] \ar[r]^-{cl} & \Omega^{j-1} \SO_{n-j}
}$$

\section{Spinning and graphing in high co-dimensions}\label{path}

This section is devoted to the concepts surrounding a new proof that
$\pi_0 \K_{n,j}$ is a group, provided $n-j>2$.  The proof is quite 
simple: show that the total-space of the fibration 
$\K_{n,j} \to \PE_{n,j} \to \K_{n-1,j-1}$ is connected.
This forces the boundary map $\pi_1 \K_{n-1,j-1} \to \pi_0 \K_{n,j}$ from
the homotopy long exact sequence to be an epi-morphism.
Showing that $\PE_{n,j}$ is connected reduces to showing
that every neat embedding of $D^j$ in $D^n$ is isotopic (through neat embeddings)
to a linear inclusion. The remainder of the section elaborates on 
ingredients used in the proof and its consequences.  The boundary map  
$\Omega \K_{n-1,j-1} \to \K_{n,j}$ is shown to be
homotopic an explicitly-defined graphing map $\gr_1 \co  \Omega\K_{n-1,j-1} \to \K_{n,j}$ 
in \fullref{lgh}.  Propositions \ref{rel-lith} and \ref{relate} demonstrate
that $\gr_1$ is a variant of Litherland's deform-spinning construction \cite{Lith}. 
Goodwillie's dissertation is invoked, showing that
$\gr_1$ is a surprisingly highly-connected map. This allows the computation of the
first non-trivial homotopy groups of $\K_{n,j}$ provided $2n-3j-3 \geq 0$.
Using some computations of Victor Turchin and a quadrisecants argument,
an explicit generator is constructed for $\pi_{2n-6} \K_{n,1}$. Via spinning, 
this gives new explicit constructions of Haefliger's spheres $\pi_0 \K_{n,j}$ for
$2n-3j-3=0$. 

The next proposition is an old result which is known to hold in far
greater generality (see Hudson \cite{Hud} and Goodwillie \cite{GoodD}).
Goodwillie's generalisation will later be used in this paper. So strictly
speaking, this proposition is redundant.  The proof is included as
several later developments in this section build on it, making it the
natural starting point.

\begin{prop}\label{Embdk}
Provided $n-j>2$, the map $\pi_1 \K_{n-1,j-1} \to \pi_0 \K_{n,j}$ is an epi-morphism.
The spaces $\Emb(D^j,D^n)$ and $\PE_{n,j}$ are connected.
\end{prop}
\begin{proof}
Once $\Emb(D^j,D^n)$ is shown to be connected, 
the remaining results follow from the homotopy long exact 
sequences of the fibrations $\K_{n,j} \to \PE_{n,j} \to \K_{n-1,j-1}$ and
$\PE_{n,j} \to \Emb(D^j,D^n) \to V_{n,j}$ from  \fullref{Splitprop}.

\begin{itemize}
\item Consider $n = 4$. The path-connectivity of $\Emb(D^1,D^4)$ is 
well-known and appears in many places. Let $f \in \Emb(D^1,D^4)$, and 
isotope it to be standard on the boundary: $f(-1)=(-1,0,0,0)$ and $f(1)=(1,0,0,0)$.  
Let $v \in S^3$. By Sard's theorem, the projection of $f$ into the orthogonal complement of 
$v$ is generically an embedding. Choose one such value for $v$ such that
$c = \langle v, (1,0,0,0)\rangle  > 0$.  Then the formula 
$f(t) - a\langle f(t),v\rangle v + act\cdot v$ describes a path
(parametrised by $a \in [0,1]$) in $\Emb(D^1,D^4)$, starting at $f$ and
ending at a function which is monotone increasing in the direction of $v$, thus isotopic
to $t \longmapsto (t,0,0,0)$ by the straight-line homotopy.
\item Consider $n=5$. As in the previous case, isotope $f \in \Emb(D^2,D^5)$ 
to be standard on the boundary, and let $f_a \co  D^2 \to D^5$ for $a \in [0,1]$
be the straight-line homotopy from $f$ to the standard inclusion.  By the 
weak Whitney immersion theorem, one can assume $f_a$ is generically an embedding,
with only finitely many times $a$ for which it has an isolated, regular double point.
Wu \cite{Wu} developed a 1--parameter `Whitney trick' for this situation, 
to remove the double points from the family. 
\item Consider the case $n \geq 6$ and let $e \co  D^j \to D^n$ be a proper 
embedding.  Let $B \subset D^j$ be the open ball of radius $\frac{1}{2}$, 
centred about the origin.  Consider $D^j = D^j \times \{0\}^{n-j} \subset D^n$.
By a local linearisation, isotope $e$ so that it agrees with inclusion on
$B$, $e(x)=x$ for all $x \in B$. Let $U$ be the open ball of radius
$\frac{1}{2}$ centred about $0$ in $D^n$, and isotope $e$ so that
$e(D^j) \cap U = e(B)$.  
Let $W = D^n \setminus U$, $W_1 = \partial U$ and $W_2 = \partial D^n$.  
$\partial W = W_1 \sqcup W_2$.  $W_i \to W$ is a homotopy-equivalence 
for $i\in \{1,2\}$, since $W$ is a product. Let $V = e(D^j \setminus B)$ with 
$V_1 = W_1 \cap V$, $V_2 = W_2 \cap V$, and let 
$f \co  V_1 \times \bigl[\frac{1}{2},1\bigr] \to W$ be the map defined by 
$f(v,t)=e(2tv)$. $f$ maps $V_1 \times \bigl[\frac{1}{2},1\bigr]$ 
diffeomorphically to $V$.  Smale \cite[Corollary 3.2]{Smale2} states that $f$ extends 
to a diffeomorphism of pairs $f \co  (W_1,V_1) \times
\bigl[\frac{1}{2},1\bigr] \to (W,V)$.
Therefore it further extends to a diffeomorphism of pairs 
$f \co  (D^n,D^j) \to (D^n,\img(e))$. So $e = f\circ h$ where $h$ is the standard inclusion 
$h\co  D^j \to D^n$.  Given an orientation-preserving diffeomorphism $f$ of $D^n$ it 
acts on $\Emb(D^j,D^n)$, but the action is trivial on $\pi_0 \Emb(D^j,D^n)$ -- the idea is that one can linearise $f$ on the complement of a neighbourhood of a point in the boundary of $D^n$ (a similar argument is given in the proof of
\fullref{nullhlem}). \proved
\end{itemize}
\end{proof}

The earliest claim in the literature that $\Emb(D^j,D^n)$ is connected
for $n-j>2$ seems to be made by Haefliger. It appears in his AMS math
review of Zeeman's paper \cite{Zeeman}. Perhaps the above
proof is similar to what Haefliger had in mind, as he states the result
follows from Smale's paper \cite{Smale2}. It would be interesting to
know if there are any more elementary proofs.

The fibre-sequence $\K_{n,j} \to \PE_{n,j} \to \K_{n-1,j-1}$
`backs-up' to a fibre-sequence
$$\Omega \K_{n-1,j-1} \to \K_{n,j} \to \PE_{n,j}$$
by \fullref{hom-fib}.  The remainder of this section is devoted
to the properties of the `connecting map' $\Omega \K_{n-1,j-1} \to \K_{n,j}$
and its relatives. 

\begin{prop}\label{lgh}
The connecting-map $\Omega \K_{n-1,j-1} \to \K_{n,j}$ is homotopic to
$$ 
\xymatrix@R=10pt{  \Omega \K_{n-1,j-1} \ar[r]^-{\gr_1}  \ar@{}[d]|-{\upepsilon} & \K_{n,j} \ar@{}[d]|-{\upepsilon} \\
f \ar@{|->}[r] & \left[(t_0,t_1,\ldots,t_{j-1}) \longmapsto 
                 \left(t_0,f(t_0)(t_1,\ldots,t_{j-1})\right)\right]
}$$
and the connecting map $\Omega \EK{j-1,M} \to \EK{j,M}$ is
homotopic to
$$ 
\xymatrix@R=10pt{  \Omega \EK{j-1,M} \ar[r]^-{\gr_1}  \ar@{}[d]|-{\upepsilon} & \EK{j,M} \ar@{}[d]|-{\upepsilon} \\
f \ar@{|->}[r] & \left[(t_0,t_1,\ldots,t_{j-1},m) \longmapsto 
\left(t_0,f(t_0)(t_1,\ldots,t_{j-1},m)\right)\right].
}$$
\end{prop}
\begin{proof}
The two cases are essentially the same, so restrict attention to the fibration
$$\xymatrix{\EK{j,M} \ar[r]^-i & \PK{j,M} \ar[r]^-p & \EK{j-1,M}}.$$ 
By \fullref{hom-fib}
$$HF(i) = \{ f \co  [0,1] \to \PK{j,M}, f(0)=\Id_{\Real^j\times M}, 
f(1) \in \EK{j,M} \}.$$ 
The map $HF(i) \to \Omega \EK{j-1,M}$ defined in \fullref{hom-fib} is a
weak homotopy equivalence.  Palais has proved that every embedding space
has the homotopy-type of a CW--complex (see Palais \cite{P2}).  Strictly
speaking, he proves embedding spaces are dominated by CW--complexes,
but at that time it was a well-known theorem of Whitehead's that a space
dominated by a CW--complex has the homotopy-type of a (perhaps different)
CW--complex \cite{Whitehead}. The further fact that the various loop space
and homotopy-fibre constructions send spaces with the homotopy-type
of CW--complexes to spaces having the homotopy-type of CW--complexes is
due to Milnor \cite{Milnor}.  Thus, $HF(i) \to \Omega \EK{j-1,M}$ is
a homotopy-equivalence.

An explicit homotopy-inverse of $\Omega \EK{j-1,M} \to HF(i)$ is exhibited.
Given $f \in \Omega \EK{j-1,M}$, consider the object
$$(t, t_1,\ldots,t_j,m) \longmapsto
\left\{ \begin{array}{ll}
        (t_1,f(t_1)(t_2,\ldots,t_j,m)) & \text{for } 2t-1 \leq t_1\\
        (t_1,f(2t-1)(t_2,\ldots,t_j,m)) & \text{for } t_1 \leq 2t-1
\end{array}
\right.$$
This would be the `right' map $\Omega \EK{j-1,M} \to HF(i)$ 
(with loop-space parameter $t$) if it was a smooth function in the variable $t_1$.  
Consider a smooth `wet blanket' function $b \co  \Real \to \Real$ with the properties:
\begin{itemize}
\item $b(x) = x$ for all $x \leq 0$
\item $b(x) = 1/2$ for all $x \geq 1$
\item $b'(x) \geq 0$ for all $x \in \Real$.
\end{itemize}
Such a function can be obtained in closed-form as
$$b(x) = \int_0^x \left( 1- \int_0^x B(x) dx \right) dx $$
where $B \co  \Real \to \Real$ is any smooth function such that
$B\bigl(\frac{1}{2}+x\bigr)=B\bigl(\frac{1}{2}-x\bigr)$ and $B(x) \geq 0$ for all $x \in \Real$,  
with $B(x)=0$ for all $|x-\frac{1}{2}|\geq \frac{1}{2}$ and 
$\int_{-\infty}^\infty B(x) dx = 1$.  

For $t \in \Real$ define $b_t \co  \Real \to \Real$ as
$b_t(x) = b(x-t)+t$. Consider the
function $\Omega \EK{j-1,M} \to HF(i)$ defined by sending
$f \in \Omega \EK{j-1,M}$ to $\tilde f \in HF(i)$ by the formula
\begin{equation}
\label{eqn:star}
\tilde f(t)(t_1,\ldots,t_j,m)=
\bigl(t_1,f\bigl(b_{\frac{-3+5t}{2}}(t_1)\bigr)(t_2,t_3,\ldots,t_j,m)\bigr)
\tag{$*$}
\end{equation}
The composite $\Omega \EK{j-1,M} \to HF(i) \to \Omega \EK{j-1,M}$
is obtained by setting $t_1=1$ in \eqref{eqn:star}, thus $f$ is sent to the map
$$\bigl[(t,t_2,\ldots,t_j,m) \longmapsto
f\bigl(b_{\frac{-3+5t}{2}}(1)\bigr)(t_2,\ldots,t_j,m)\bigr] \in \Omega \EK{j-1,M}$$
which is just a reparametrisation of $f$ by $b_{\frac{-3+5t}{2}}(1)$ (thought of
as a function of $t$). Since $b_{\frac{-3+5t}{2}}(1)$ is an increasing function of
$t$ it is homotopic to the identity.
\end{proof}

Zeeman proved that the complements of certain co-dimension two `twist-spun' 
knots fibre over $S^1$ \cite{ZeemTwist}.  Litherland later went on to formulate 
a more general notion of spinning, at the time called `deform-spinning,' further 
generalising Zeeman's theorem to this context \cite{Lith}. The Zeeman--Litherland
results are important for a number of reasons -- one being that they are an excellent
source of embeddings of $3$--manifolds in $S^4$, as the Seifert-surfaces of embeddings
of $S^2$ in $S^4$. The next proposition points out that the connecting map 
$\gr_1 \co  \Omega \K_{n-1,j-1} \to \K_{n,j}$ is a mild variation of Litherland's 
spinning construction. 

\begin{figure}\label{fig2}
\centering
\labellist\small
\pinlabel {$P_n$} [b] at 720 333
\endlabellist
\includegraphics[width=10cm]{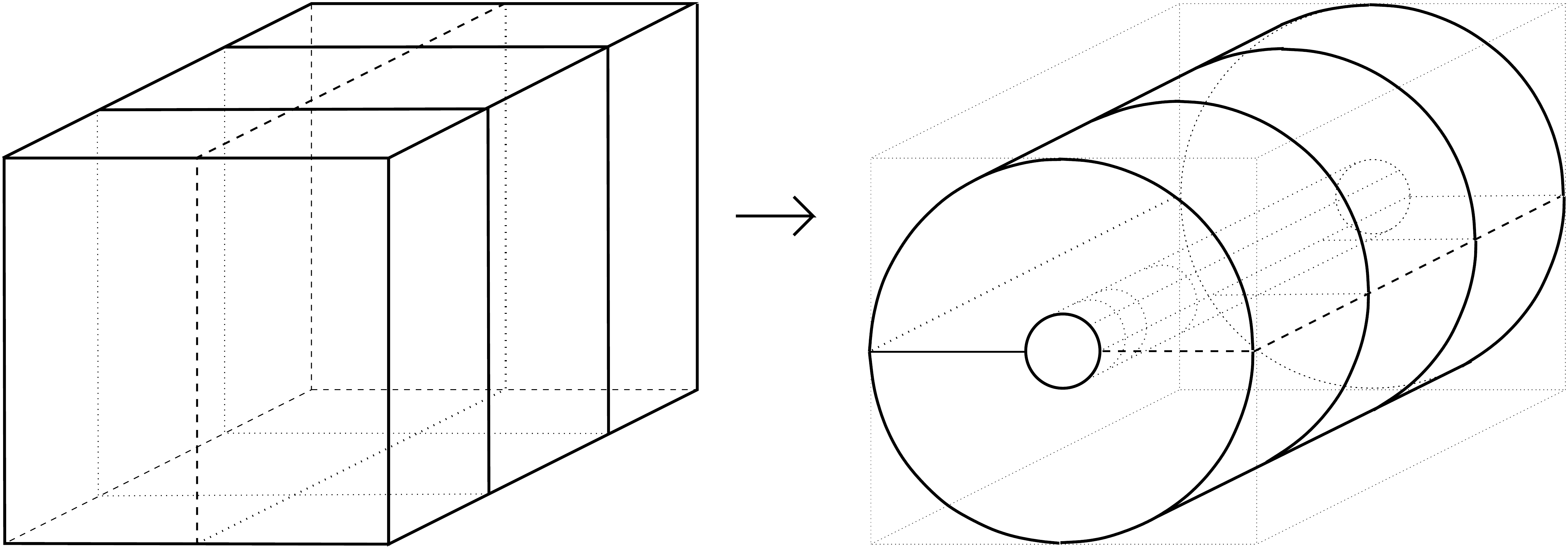}
\caption{}
\end{figure}

Given a topological space $X$, denote the space of
continuous functions $f \co  S^1 \equiv \Real/2\Zed \to X$ by $LX$
called the `free loop space of $X$.'
Define $P_2\co  \I^2 \to \I^2$ by
$P_2(t_1,t_2)=\bigl(\frac{t_2+2}{3} \cos(\pi t_1),\frac{t_2+2}{3} \sin(\pi
t_1)\bigr)$ 
and $P_n \co  \I^n \to \I^n$ as $P_n = P_2 \times \Id_{\I^{n-2}}$. See
\fullref{fig2}.
Notice $P_n$ is an embedding on the interior of $\I^n$, and is globally one-to-one 
except for the equality $P_n(-1,t_2,t_3,\ldots, t_n)=P_n(1,t_2,\ldots,t_n)$.
\begin{defn}\label{litherland-spin}
Given $f \in L\K_{n-1,j-1}$, let $h \co  \Real^j \to \Real^n$
be the function
$$h(t_0,t_1,\ldots,t_{j-1}) = (t_0,f(t_0)(t_1,\ldots,t_{j-1})),$$
and consider the composite $P_n \circ h \circ P_j^{-1}$.
It is well-defined on the image of $P_j$. On $\partial P_j(\I^j)$ 
it agrees with the standard inclusion 
$\Real^j \to \Real^n$. Define $\gr_1(f) \in \K_{n,j}$ to be the unique extension of 
$P_n \circ h \circ P_j^{-1}$ such that $\gr_1(f)_{|\Real^j \setminus P_j(\I^j)}$ 
agrees with the standard inclusion.
\end{defn}

\begin{prop}\label{rel-lith}
The diagram
$$\xymatrix{
L\K_{n-1,j-1} \ar[r]^-{\gr_1} & \K_{n,j} \\
\Omega \K_{n-1,j-1} \ar[u] \ar[ur]_{\gr_1} & } $$
is homotopy-commutative.
\end{prop}
\begin{proof}
There exists a $1$--parameter family $P_n(t) \co  \I^n \to \I^n$ 
for $t \in [0,1]$ satisfying
$P_n(0) = P_n$, $P_n(1) = \Id_{\I^n}$, such that for all
$t \in (0,1]$ the function $P_n(t) \co  \I^n \to \I^n$ is an embedding.  
Substituting $P_n(t)$ for $P_n$ in the definition of 
$\gr_1 \co  L \K_{n-1,j-1} \to \K_{n,j}$ gives the desired homotopy.
\end{proof}

In the literature, Litherland spinning is not defined as the map
$\gr_1 \co  L\K_{n-1,j-1} \to \K_{n,j}$, but what Litherland defined in
\cite{Lith}, when appropriately adapted to the smooth category,
turns out to be precisely $\gr_1$. This is the content of
\fullref{relate}. 

$\EK{n,*}$ is the group of diffeomorphisms of $\Real^n$
whose support is contained in $\I^n$, thus it acts (by composition on the
left) on $\K_{n,j}$. Notice that if $n-j>0$, $f \in \K_{n,j}$ and $g \in \EK{n,*}$
then $g \circ f$ is in the the same path-component of $\K_{n,j}$ as $f$. 
In fact, much more is true. Let $\K_{n,j}(f)$ denote the path-component of
$f$ in $\K_{n,j}$.

\begin{lem}\label{nullhlem}
Provided $n-j>0$ and $f \in \K_{n,j}$, the map 
$\EK{n,*} \to \K_{n,j}$ given by sending $g \in \EK{n,*}$
to $g \circ f$ is a null-homotopic fibration whose image is
$\K_{n,j}(f)$. The fibre of this fibration is denoted $\Diff(\I^n,f)$.
\end{lem}
\begin{proof}
That the map is a fibration is classical (see Cerf \cite{Cerf}).  That the image
contains $\K_{n,j}(f)$ follows from the isotopy extension theorem. 
Consider an orientation-preserving affine-linear transformation
$L \co  \Real^n \to \Real^n$ such that $L(\I^n) \subset \I^n$.  Given
$g \in \EK{n,*}$ notice that $L \circ g \circ L^{-1} \in \EK{n,*}$,
moreover the support of $L \circ g \circ L^{-1}$ is contained in
$L(\I^n)$. The space of orientation-preserving affine linear transformations
of $\Real^n$ which preserves $\I^n$ is connected, thus there is a path $L_t$
in this space such that $L_0 = \Id_{\Real^n}$ and $L_1 = L$.  The function
$$ 
\xymatrix@R=10pt{  [0,1] \times \EK{n,*} \ar[r]  \ar@{}[d]|-{\upepsilon} & 
                  \K_{n,j} \ar@{}[d]|-{\upepsilon} \\
(t,g) \ar@{|->}[r] &  L_t \circ g \circ L_t^{-1} \circ f }
$$
is a null-homotopy of the map $\EK{n,*} \to \K_{n,j}$ 
provided $L(\I^n) \cap f(\Real^j) = \emptyset$, which can always be arranged provided
$n-j>0$, by Sard's theorem.
\end{proof}

The map $\pi_1 \K_{n,j}(f) \to \pi_0 \Diff(\I^n,f)$ is therefore a bijection
onto the subgroup of $\pi_0 \Diff(\I^n,f)$ which is the kernel of the forgetful
map $\pi_0 \Diff(\I^n,f) \to \pi_0 \EK{n,*}$. 
Given an element $g \in \pi_1 \K_{n,j}(f)$, let $\tilde g \in \pi_0 \Diff(\I^n,f)$ 
be its image.  Given $g \in \pi_1 \K_{n,j}(f)$ and $\gr_1 g \in \K_{n+1,j+1}$
denote the one-point compactification by $\overline{\gr_1 g} \in \Emb(S^{j+1},S^{n+1})$.

Starting from an element $h \in \Diff(\I^n,f)$ which is in the kernel of the forgetful map
$\Diff(\I^n,f) \to \pi_0 \EK{n,*}$, Litherland gave a `surgery' description
\cite{Lith} of an embedding $S^{j+1} \to S^{n+1}$.  Consider $\I^{n+2}$ to be the 
product $\I^{n+2} = \I^n \times \I^2$, so $\partial \I^{n+2} = 
\I^n \times (\partial \I^2) \cup (\partial \I^n) \times \I^2$. 
Think of $\I^n \times (\partial \I^2)$ as a trivial $\I^n$--bundle over 
$\partial \I^2$, therefore it is diffeomorphic to the bundle over 
$\partial \I^2$ with fibre $\I^n$ and monodromy given
by $h$. Call this space $\I^n \times_h \partial I^2$.
Since $h$ acts as the identity on $\partial \I^n$, the boundary of  
$\I^n \times_h \partial \I^2$ is canonically identified with 
$\partial \I^n \times \partial \I^2$.  Thus the union 
$$\bigl((\I^n,f) \times_h \partial \I^2\bigl) \cup (\partial \I^n, \partial \I^j)\times \I^2$$
makes sense as a manifold pair.  Identify $\partial \I^{n+2}$ with
$S^{n+1} \subset \Real^{n+2}$ by radial projection from the origin. 
Thus, $\left((\I^n,f) \times_h \partial \I^2\right) \cup (\partial \I^n, \partial \I^j)\times \I^2$
describes an embedding of $S^{j+1}$ in $S^{n+1}$. This is Litherland's
deform-spun knot construction \cite{Lith}.

\begin{prop}\label{relate}
Given $g \in \pi_1 \K_{n,j}(f)$, the `Litherland spun' knot 
$\left((\I^n,f) \times_{\tilde g} \partial \I^2 \right) \cup 
 (\partial \I^n, \partial \I^j)\times \I^2$
and $\overline{\gr_1 g} \in \Emb(S^{j+1},S^{n+1})$ are isotopic, once
$S^{n+1}$ is identified with $\partial \I^{n+2}$ via radial projection.
\end{prop}
\begin{proof}
The key step is to remember that the identification of 
$\I^n \times (\partial \I^2)$ with
$\I^n \times_{\tilde g} \partial I^2$ is made via the null-isotopy of
$\tilde g$ when considered as an element of $\EK{n,*}$. Under this identification,
the two definitions are identical. 
\end{proof}

Given $f \in \K_{n,j}$ and $g \in \Omega \K_{n,j}(f)$, let $C_f$ be the complement
of an open tubular neighbourhood of $\wbar{f}$ in $S^n$. 
By the above argument, the complement of $\overline{\gr_1(g)}$ in 
$S^{n+1}$ is diffeomorphic to $C_f \rtimes_{\tilde g} S^1$ union 
a $2$--handle and an $(n{-}j{+}1)$--handle. 
Here $C_f \rtimes_{\tilde g} S^1$ indicates the $C_f$ bundle over
$S^1$ with monodromy induced by $\tilde g$. This gives a presentation
$$\pi_1 C_{\gr_1(g)} \simeq \pi_1 C_f / \langle {\tilde g}.x=x \ \forall x \in \pi_1 C_f \rangle$$
where $\langle {\tilde g}.x=x \ \forall x \in \pi_1 C_f \rangle$ is the normal subgroup of
$\pi_1 C_f$ generated by the relations ${\tilde g}.x=x$ for all $x \in \pi_1 C_f$. 

\begin{eg}
If $g \in \Omega\K_{3,1}(f)$ is the Gramain element (rotation by $2\pi$ about 
the long axis), its action on $\pi_1 C_f$ is conjugation by the meridian.
Thus $\pi_1 C_{\gr_1(g)}$ is trivial, as all knot groups are `normally
generated' by a meridian.  This observation anticipates the Zeeman--Litherland
theorem, which states that $\gr_1(g)$ is the unknot (see Zeeman
\cite{ZeemTwist} and Litherland \cite{Lith} )
whenever $g$ is the Gramain element.  The Zeeman--Litherland theorem is
stated in full generality in \fullref{survey}.
\end{eg}

The spaces $\K_{n,n} = \EK{n,*}$ are the groups of diffeomorphisms of a cube,
and have the homotopy-type of $\Diff(D^n)$, the group of diffeomorphisms
of a disc which are the identity on the boundary.
The maps $\gr_1 \co  \Omega \K_{n,n} \to \K_{n+1,n+1}$ have been studied in this
context. Define $\gr_2 \co  \Omega^2 \K_{n,j} \to \K_{n+2,j+2}$ to be
the composite $\gr_1 \circ \Omega \gr_1$ where 
$\Omega \gr_1 \co  \Omega^2 \K_{n,j} \to \Omega \K_{n+1,j+1}$ is the induced map
of $\gr_1$. Similarly define $\gr_i \co  \Omega^i \K_{n,j} \to \K_{n+i,j+i}$. 
In the literature (see Antonelli, Burghelea and Kahn \cite{ABK2}, Weiss
and Williams \cite{Weiss} and Gromoll \cite{Grom}) elements of $\pi_0
\K_{n,n}$ which are in the image of $\gr_i \co  \pi_i \K_{n-i,n-i}
\to \pi_0 \K_{n,n}$ but which are not in the image of $\gr_{i+1}$ are
typically said to have Gromoll degree $i$.

\begin{defn}
An element $f \in \pi_0 \K_{n,j}$ has (Gromoll) degree
$i$ if it is in the image of the $i$th graphing map 
$\gr_i \co  \pi_i \K_{n-i,j-i} \to \pi_0 \K_{n,j}$ but not
in the image of the $(i{+}1)$st graphing map $\gr_{i+1}$.
\end{defn}

\begin{prop}\label{Connectivityknj}
\begin{enumerate}
\item The Gromoll degree of the elements of $\pi_0 \K_{n,j}$
is at least $2n-2j-4$ for all $n \geq j > 0$. 
\item $\K_{n,j}$ is $(2n{-}3j{-}4)$--connected for all $n \geq j \geq 1$.
Provided $2n-3j-3 \geq 0$ and $n-j>2$ the
first non-trivial homotopy group of $\K_{n,j}$ is
$$
\pi_{2n-3j-3} \K_{n,j} \simeq
\left\{
\begin{array}{ll}
\Zed & j=1 \text{ or } n-j \text{ is odd } \\
\Zed_2 & j > 1 \text{ and } n-j \text{ is even }
\end{array}
\right.
$$
The elements of $\pi_0 \K_{n,j}$ for $2n-3j-3=0$ have Gromoll degree
$(j-1)$, ie: $\gr_{j-1} \co  \pi_{j-1} \K_{n-j+1,1} \to \pi_0 \K_{n,j}$ is onto.
\item $\Emb(S^j,S^n)$ is $\min\{(2n-3j-4),(n-j-2)\}$--connected 
for all $n \geq j \geq 1$. Let $m=\min\{2n-3j-3,n-j-1\}$. Provided 
$2n-3j-3\geq 0$ and $n-j>2$ the first non-trivial homotopy-group of $\Emb(S^j,S^n)$ is
$$
\pi_m \Emb(S^j,S^n) \simeq
\left\{
\begin{array}{ll}
\Zed & 2n-3j-3 < n-j-1, (j=1 \text{ or } n-j \text{ odd}) \\
\Zed & 2n-3j-3 > n-j-1, n-j \text{ even } \\
\Zed_2 & 2n-3j-3 < n-j-1, j>1 \text{ and } n-j \text{ even} \\
\Zed_2 & 2n-3j-3 > n-j-1, n-j \text{ odd }        \\
\Zed \oplus \Zed_2 & 2n-3j-3=n-j-1 
\end{array}
\right.
$$
\item $\Emb(S^j,\Real^n)$ is $\min\{2n-3j-4,n-j-2\}$ connected
for all $n \geq j + 2 \geq 3$. Let $m=\min\{2n-3j-3,n-j-1\}$.
Provided $2n-3j-3\geq 0$ and $n-j>2$ the first non-trivial homotopy group of
$\Emb(S^j,\Real^n)$ is 
$$
\pi_m \Emb(S^j,\Real^n) \simeq
\left\{
\begin{array}{ll}
\Zed & 2n-3j-3 < n-j-1, (j=1 \text{ or } n-j \text{ odd}) \\
\Zed_2 & 2n-3j-3 < n-j-1, j>1 \text{ and } n-j \text{ even} \\
\Zed & 2n-3j-3 > n-j-1  \\
\Zed^2 & 2n-3j-3 = n-j-1, (j=1 \text{ or } n-j \text{ odd}) \\
\Zed \oplus \Zed_2 & 2n-3j-3 = n-j-1, j>1 \text{ and } j \text{ even} \\
\end{array}
\right.
$$
\item $\PE_{n,j}$ is $(2n{-}2j{-}5)$--connected for all $n-j>2$.
\item $\Emb(D^j,D^n)$ is $(n{-}j{-}2)$--connected for all $n - j > 2$.
\end{enumerate}
\end{prop}
\begin{proof}
(5)\qua That $\PE_{n,j}$ is $2n-3j-5$ connected follows directly from Goodwillie's 
dissertation \cite[Theorem~C, page~9]{GoodD}. 

(6)\qua This result follows from (5) and \fullref{Splitprop}.

(1)\qua Consider the homotopy fibre-sequence 
$\Omega \K_{n-1,j-1} \to \K_{n,j} \to \PE_{n,j}$ from \fullref{lgh}.
Since $\PE_{n,j}$ is $(2n{-}2j{-}5)$--connected, $\pi_1 \K_{n-1,j-1} \to \pi_0 \K_{n,j}$
is epic for $n-j>2$.  Moreover, $\pi_2 \K_{n-2,j-2} \to \pi_1 \K_{n-1,j-1}$
is also epic, as $\pi_1 \PE_{n-1,j-1}$ is trivial.  The result follows by induction.

(2)\qua There is a computation of the 3rd stage of the Goodwillie
tower for $\K_{n,1}$ in \cite{Bud}. This is a $(2n{-}6)$--connected
map $\K_{n,1} \to AM_3$. $AM_3$ is known to have the
homotopy-type of the 3--fold loop-space on the homotopy fibre of the inclusion
$S^{n-1} \vee S^{n-1} \to S^{n-1} \times S^{n-1}$, thus 
$\K_{n,1}$ is $(2n{-}7)$--connected. The first non-trivial 
integral homology group of $\K_{n,1}$ is computed by Victor Turchin \cite{Tur2}
(see the computations for the homology of the complexes 
$CT_0D^\text{even}$ and $CT_0D^\text{odd}$ for $j=4$, $i=2$). 
Turchin's result is that $H_{2n-6}(\K_{n,1};\Zed) \simeq \Zed$,
so by the Hurewicz Theorem, $\pi_{2n-6} \K_{n,1} \simeq \Zed$. 
That verifies the result for $\K_{n,1}$.  

Consider the space
$\K_{n+j,j+1}$ for $j \geq 1$. The fibre-sequence
$$\Omega \K_{n+j-1,j} \to \K_{n+j,j+1} \to \PE_{n+j,j+1}$$
has a $(2n{-}7)$--connected base-space. In the special
case of $j=1$ the fibre 
has first non-trivial homotopy group in dimension $2n-7$. But
$\pi_{2n-7} \PE_{n+1,2}$ is trivial, thus 
$\pi_{2n-6}\K_{n,1} \to \pi_{2n-7} \K_{n+1,2}$
is epic with kernel generated by the image of $\pi_{2n-6} \PE_{n+1,2}$, 
giving the isomorphism
$$\pi_{2n-7} \K_{n+1,2} \simeq \pi_{2n-6}\K_{n,1} / \img\left(\pi_{2n-6}
\PE_{n+1,2}\right).$$
Repeat the argument for $j>1$, inductively assuming that
the first non-trivial homotopy group of $\Omega \K_{n+j-1,j}$
is $\pi_{2n-j-6} \Omega \K_{n+j-1,j}$ and isomorphic to 
$$\pi_{2n-6}\K_{n,1} / \img\left(\pi_{2n-6} \PE_{n+1,2}\right).$$ 
Since $\PE_{n+j,j+1}$ is $(2n{-}7)$--connected, the map
$$\pi_{2n-j-6} \Omega \K_{n+j-1,j} \to \pi_{2n-j-6} \K_{n+j,j+1}$$
is an isomorphism of first non-trivial homotopy-groups, thus
for all $j \geq 1$ there is an isomorphism 
$$\pi_{2n-j-6} \K_{n+j,j+1} \simeq \pi_{2n-6}\K_{n,1} /
\img\left(\pi_{2n-6} \PE_{n+1,2}\right).$$
Setting $j$ equal to $2n-6$ gives the isomorphism
$$\pi_0 \K_{3n-6,2n-5} \simeq \pi_{2n-6}\K_{n,1} / \img\left(\pi_{2n-6} \PE_{n+1,2}\right).$$
Haefliger's computations \cite{Haefliger2} completes the proof:
$$\pi_0 \K_{3n-6,2n-5} \simeq \left\{
\begin{array}{ll}
\Zed_2 & \text{for } n \geq 4 \text{ odd} \\
\Zed & \text{for } n \geq 4 \text{ even.}
\end{array} \right.
$$ 
(3)\qua \fullref{Splitprop} gives us a homotopy equivalence
$\Emb(S^j,S^n) \simeq \SO_{n+1} \times_{\SO_{n-j}} \K_{n,j}$.
Since $\SO_{n+1}/\SO_{n-j} \equiv V_{n+1,j+1}$ is
$(n{-}j{-}1)$--connected, the homotopy long exact sequence of the fibration
$\K_{n,j} \to \Emb(S^j,S^n) \to V_{n+1,j+1}$ tells us that
$\Emb(S^j,S^n)$ is $\min\{n{-}j{-}1,2n{-}3j{-}4\}$--connected.
Since  the bundle
$$\Emb(S^j,S^n) \to V_{n+1,j+1}$$ is
split, the first non-trivial homotopy group of
$\Emb(S^j,S^n)$ can be computed directly.

(4)\qua For $\Emb(S^j,\Real^n)$ use the homotopy equivalence
$\Emb(S^j,\Real^n) \simeq \SO_n \times_{\SO_{n-j}} 
\left( C \rtimes \K_{n,j} \right)$ from \fullref{long-euclid}.  
The bundles 
$C \rtimes \K_{n,j} \to \K_{n,j}$ and
$\SO_n \times_{\SO_{n-j}} \left( C \rtimes \K_{n,j} \right) \to
V_{n,j}$ are split, so the computation follows directly.
\end{proof}

An interesting corollary is that there are `exotic families' of
smooth $2$--discs in the $6$--disc.

\begin{cor}
$\pi_{2n-6} \PE_{n+1,2}$ has rank at least $1$ provided $n \geq 5$
is odd. 
\end{cor}

Brian Munson \cite{Munson} gave a lower bound of $\min\{2n-3j-4,n-j-2\}$
on the connectivity of $\Emb(S^j,\Real^n)$.  \fullref{Connectivityknj}
proves that Munson's lower bound is sharp.

\begin{figure}\label{fig3} \centering
\labellist\small \pinlabel {$\Real^n$} [t] at 330 300 \pinlabel
{$P_1$} [b] at 150 238 \pinlabel {$P_2$} [b] at 525 223 \endlabellist
\includegraphics[width=10cm]{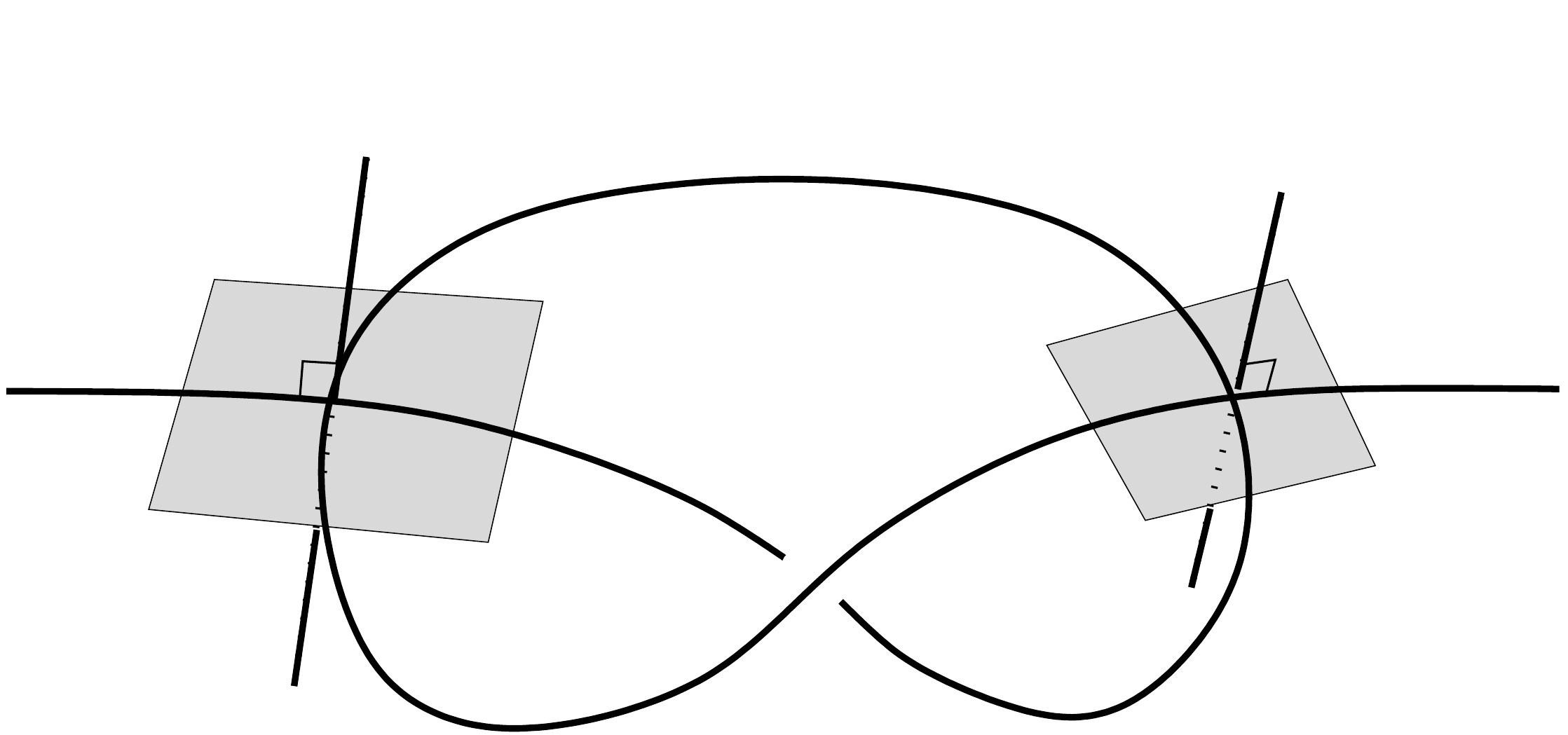} \caption{} \end{figure}

The rest of this section is devoted to a geometric construction of
the generators of $\pi_{2n-6} \K_{n,1}$ for $n \geq 4$.  Take a `long'
immersion $f\co \Real \to \Real^3 \subset \Real^n$ having two regular
double points $f(t_1)=f(t_3)$, $f(t_2)=f(t_4)$ with $t_1 < t_2 < t_3 <
t_4 \in \Real$ such that one of the four resolutions of $f$ in $\Real^3$
is a trefoil knot.  Let $Tf_i$ be the tangent space to $f(\Real)$ at
$t_i$.  Let $P_1$ be the orthogonal complement to $Tf_1 \oplus Tf_3$
in $\Real^n$, and $P_2$ the orthogonal complement of $Tf_2 \oplus
Tf_4$ in $\Real^n$.  $P_1$ and $P_2$ are $(n{-}2)$--dimensional,
so if $S_1$ and $S_2$ are the unit sphere of $P_1$ and $P_2$
respectively they are both $(n{-}3)$--dimensional. There is a
`resolution function' $r \co  S_1 \times S_2 \to \K_{n,1}$ given
by perturbing $f$ near the double points via bump-functions whose
directions are prescribed by the pair $(v_1,v_2)\in S_1 \times S_2$.
See \fullref{fig3}.  The claim is that $r$ is a generator of $H_{2n-6}
(\K_{n,1};\Zed) \simeq \Zed$.  
One could potentially trace through the computations of Turchin
\cite{Tur2} and Vassiliev \cite{Vass2} to verify that $r$ generates
$H_{2n-6} (\K_{n,1};\Zed)$.  The following approach is perhaps more
direct.  It is inspired by the author's quadrisecant description of the
type-2 Vassiliev invariant for knots $\Real^3$ \cite{Bud}. The idea
is to construct an integral co-homology class $\nu_2 \in H^{2n-6}
(\K_{n,1};\Zed)$ such that if $x \in H_{2n-6} (\K_{n,1};\Zed)$ is
represented as an oriented $(2n{-}6)$--dimensional manifold mapping into
$\K_{n,1}$ then $\nu_2(x)$ can be computed as a signed count of the number
of alternating quadrisecants along the family of long knots represented
by $x$.  Every class in $H_{2n-6}(\K_{n,1};\Zed)$ is realisable as a map
from an oriented $(2n{-}6)$--dimensional manifold $M$ to $\K_{n,1}$ since
$\K_{n,1}$ is $(2n{-}7)$--connected (\fullref{Connectivityknj}). Moreover,
by the Hurewicz theorem, $M$ can be assumed to be $S^{2n-6}$, as
$\pi_{2n-6}\K_{n,1} \simeq H_{2n-6}(\K_{n,1};\Zed)$.

\begin{defn}
Given two points $x,y \in \Real^n$ let $[x,y]$ denote
the oriented line segment in $\Real^n$, starting at 
$x$ and ending at $y$. 
 An alternating quadrisecant in $C_4(\Real^n)$ is
a point $(x_1,x_2,x_3,x_4) \in C_4(\Real^n)$
such that $[x_1,x_4] \subset [x_3,x_2]$ as an
oriented subinterval. $C_k M$ denotes the configuration space
of distinct $k$--tuples of points in $M$, 
$C_k M = \{ x \in M^k : x_i \neq x_j \ \forall \ i \neq j \}$.
Provided $M$ is a manifold, let $C_k[M]$ denotes the (real oriented) 
Fulton--Macpherson compactification of $C_k M$, as in \cite{Bud}.
$C_k[M]$ is a compact manifold, provided $M$ is compact. 
The `real oriented' Fulton--Macpherson compactification
has the property that the inclusion $C_kM \to C_k[M]$ is a 
homotopy-equivalence. 

Let $AQ_n \subset C_4[\Real^n]$ denote the closure of the set of
all alternating quadrisecants in $C_4(\Real^n)$. Let 
$C'_4[\Real] = \overline{\{t \in C_4(\Real) : t_1<t_2<t_3<t_4 \}}$.
Given $f \in \K_{n,1}$ let $AQ_n(f) \subset C'_4[\Real]$ 
denote the pull-back of $AQ_n$. More generally, if
$f \co  M \to \K_{n,1}$ is smooth, define 
$AQ_n(f) \subset M \times C_4'[\Real]$ as the pull-back of $AQ_n$.
\end{defn}

Given a closed, oriented $(2n{-}6)$--dimensional manifold $M$ and a map 
$f \co  M \to \K_{n,1}$ such that $f_* \co  M \times C'_4[\Real] \to C_4[\Real^n]$ 
is transverse to $AQ_n$, $AQ_n(f) \subset M \times C_4'[\Real]$ 
is a $0$--dimensional submanifold whose normal bundle is oriented by
the map. A well-defined integer invariant $\nu_2(f) \in \Zed$ is
defined as the signed count (of the relative orientations) of the points
in $AQ_n(f)$. The sign of each point of $AQ_n(f)$ could be computed by a formula
analogous to the one in \cite[Proposition~6.2]{Bud}. \fullref{welldef} is
the key technical lemma needed to show that $\nu_2(f)$ is an invariant
of the homology class of $f$.

Given $f \in \K_{n,1}$ let $\Gamma(f) \in (0,\infty]$ be the `cut radius'
of $f$ in $\Real^n$, defined as the supremum over all $R$ such that the
exponential map from $f$'s radius-$R$ normal disc bundle to $\Real^n$ is an 
embedding.  $\Gamma \co  \K_{n,1} \to (0,\infty]$ can be shown to be a 
continuous function, as $\Gamma(f)$ is the minimum of two continuous quantities
1) the focal radius of $f$ (which can be computed in terms of the 2nd fundamental
form of $f$) and 2) the minimum of the distances $L$ such that there 
exists two geodesics segments, each of length $L$, emanating from a point in 
$\Real^n$ and terminating in $f(\Real)$, orthogonal to the tangent space 
of $f(\Real)$. This kind of continuity argument is standard in differential
geometry, see Sakai \cite[Proposition~III.4.1]{tSakai} for example.

\begin{lem}\label{welldef}
Every $x \in H_{2n-6}(\K_{n,1};\Zed)$ represented by a manifold
$f \co  M \to \K_{n,1}$ can be perturbed so that $f_*$
is transverse to $AQ_n$. 
\end{lem}
\begin{proof}
Let $R$ be the cut radius of $f$, $R = \min\{ \Gamma(f(x)) : x \in M \}$.
Let $b \co  \Real \to \Real$ be a $C^\infty$--smooth function satisfying: 
\begin{itemize}
\item $b(x)=0$ for all $|x|\geq 1$
\item $b(x)=b(-x)$ for all $x \in \Real$
\item $\int_{-\infty}^\infty b(x)dx = 1$
\item $b'(x) > 0$ for all $-1 < x < 0$.
\end{itemize}
For $\epsilon > 0$ and $t \in \Real$ let $b_{\epsilon,t} \co  \Real \to \Real$ be defined as
$b_{\epsilon,t}(x) = \frac{1}{\epsilon}b\bigl(\frac{x-t}{\epsilon}\bigr)$.  
By a compactness argument, there exists an $m \in \Zed$ (perhaps very large) such that 
if $I_1,\ldots, I_m$ is the partition of $\I$ into $m$ equal-length sub-intervals, then 
for all $x \in M$ and $j \in \{1,2,\ldots,m\}$, $f(x)(I_j)$ is contained in the radius $R/2$ 
tubular neighbourhood of $f(x)$. 

Consider the function $\tilde f$ defined as
$$ 
\xymatrix@R=10pt{  M \times (\Real^n)^m \times \Real \ar[r]^-{\tilde f}  \ar@{}[d]|-{\upepsilon} &  \Real^n \ar@{}[d]|-{\upepsilon} \\
(x,v_1,\ldots,v_m,t) \ar@{|->}[r] &  
 f(t) + \sum_{j=1}^m b_{\frac{3}{2m},p_j}(t)v_j }
$$
where $p_j \in I_j$ is the mid-point of the interval $I_j$.
Since embeddings are an open subset of the space of all `long' 
smooth maps from $\Real$ to $\Real^n$ (see Hirsch \cite{Hirsch}), in some 
neighbourhood $U$ of $0$ in $(\Real^n)^m$, a restriction of $\tilde f$ 
can be thought of as a map $\wbar{f} \co  M \times U \to \K_{n,1}$. 
Consider a point $(x,y,t_1,t_2,t_3,t_4)$ of 
$AQ_n(\wbar{f}) \subset M \times U \times C_4'[\Real]$.
For each $i$, $t_i$ and $t_{i+1}$ cannot both be elements of some common $I_j$
since $(f(t_1),f(t_2),f(t_3),f(t_4))$ is an alternating quadrisecant. 
Thus 
$\wbar{f}_* \co  M \times U \times C'_4[\Real] \to C_4[\Real^n]$
is transverse to $AQ_n$. By the Transversality Theorem
(see Guillemin and Pollack \cite{Pollack}), $f$ can be approximated by
a map $M \to \K_{n,1}$ such that the induced map $M \times C_4'[\Real]
\to C_4[\Real^n]$ is transverse to $AQ_n$.
\end{proof}

\begin{thm}\label{resgen}
$\nu_2 \in H^{2n-6} (\K_{n,1};\Zed)$ is a well-defined cohomology class.
Moreover, $\nu_2(r) = \pm 1$, forcing $r$ to be a generator of 
$H_{2n-6} (\K_{n,1};\Zed) \simeq \Zed$.
\end{thm}
\begin{proof}
An alternating quadrisecant can never appear on $\partial (M \times
C'_4[\Real])$ nor can a $1$--parameter family of alternating quadrisecants
run off to infinity, thus, by the Transversality Extension Theorem
(see for example \cite[Chapter~2]{Pollack}) $\nu_2(f)$ is well-defined
integer invariant of the homology class of $f$.

In the picture of the `immersed trefoil' $f \co  \Real \to \Real^3
\subset \Real^n$ there are no quadrisecants, except the `degenerate'
quadrisecant that consisting of the secant between the two pairs
of double-points. Consider all the possible resolutions $r$ of this
immersed trefoil.  $r$ only has $4$ resolutions in $\Real^3 \subset
\Real^n$, so these are the only $4$ resolutions that could possibly
have quadrisecants. Moreover, only the resolution which is a trefoil in
$\Real^3$ has a quadrisecant.
\end{proof}

Since $\K_{n,1}$ is $(2n{-}7)$--connected, by the Hurewicz
Theorem $\pi_{2n-6} \K_{n,1} \simeq \Zed$ is generated by any map
$\tilde r \co  S^{2n-6} \to \K_{n,1}$ homologous to $r$.
One can explicitly construct such a map -- attachment of an
$(n{-}3)$--handle to $S_1 \times S_2 \times [0,1]$ along $S_1 \times \{*\} \times \{1\}$
gives a cobordism between $S_1 \times S_2$ and $S^{2n-6}$. 
$r_{|S_1 \times \{*\}}$ is null so $r$ extends over the cobordism.
$\tilde r$ can be chosen to be the restriction of this cobordism to $S^{2n-6}$.


\section{Actions of operads of little cubes on embedding spaces}\label{actions}

This section is devoted to the study of the iterated loop-space structures
on the embedding spaces $\K_{n,j}$ and $\EK{j,D^n}$, especially focusing
on the compatibility of these structures with Litherland spinning $\gr_1$.
The context of these results comes from the work of Boardman and Vogt
\cite{BV} and May \cite{May1,May2}. They give a very simple criterion
for recognising if a space $X$ has the homotopy-type of an $n$--fold
loop-space, being that $X$ admits an action of the operad of little
$n$--cubes, and that the induced monoid structure on $\pi_0 X$ is that of
a group.  A useful reference for operads relevant to topology, including
operads of cubes, is the book of Markl, Shnider and Stasheff \cite{MSS}.

There is an action of the operad of $j$--cubes on the spaces $\EK{j,M}$ and 
$\K_{n,j}$ given by concatenation (see \fullref{obvious}). 
The first instance of an action of the operad of $(j{+}1)$--cubes on
any space of the form $\EK{j,M}$ was given by Morlet \cite{Mor}. The Cerf--Morlet 
`Comparison Theorem' states that $\EK{j,*} \simeq \Omega^{j+1} (PL_j/O_j)$
(see Burghelea and Lashof \cite{BL} or Kirby and Siebenmann \cite{KS}
for a proof). Here $PL_j$ is the group of PL-automorphisms of $\Real^j$
(given a suitable topology) and $O_j$ is the group of linear
isometries of $\Real^j$. 
\begin{figure}\label{fig4}
\centering
\labellist\small
\pinlabel {$f\#g$} [t] at 690 450
\pinlabel {$g\#f$} [b] at 695 210
\pinlabel {$f$} [t] at 185 395
\pinlabel {$g$} [t] at 180 240
\endlabellist
\includegraphics[width=10cm]{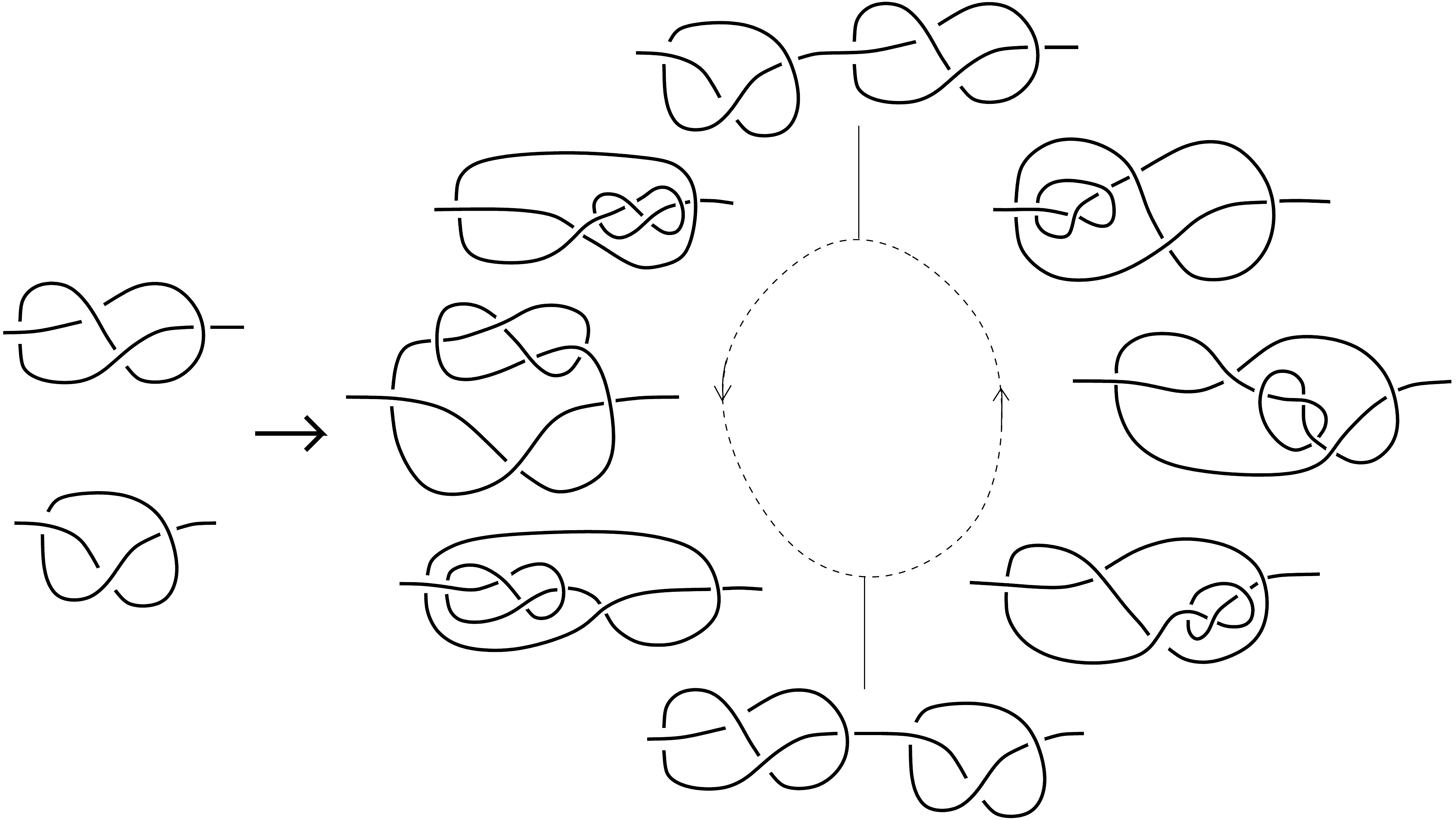}
\caption{}
\end{figure}
The first `hint' of a higher cubes action on the spaces $\EK{j,M}$
for $M$ non-trivial would perhaps be the work of Schubert \cite{Sch}.
Schubert demonstrated that the connect-sum pairing turns $\pi_0 \K_{3,1}$
into a free commutative monoid on the isotopy-classes of prime long
knots, where the demonstration of commutativity involved `pulling one
knot through another' as in \fullref{fig4}.

In `Little cubes and long knots' \cite{Cube} this idea was extended
to construct a $(j{+}1)$--cubes action on the spaces $\EK{j,M}$ for an
arbitrary compact manifold $M$.  By some elementary considerations, this
also gives an action of the operad of $(j{+}1)$--cubes on $\K_{n,j}$
for all $n-j \leq 2$. Schubert's theorem that $\pi_0 \K_{3,1}$ is a
free commutative monoid over the isotopy classes of prime long knots
generalises in this context to say that $\K_{3,1}$ is a free $2$--cubes
object over the based space $\Prime \sqcup \{*\}$ where $\Prime \subset
\K_{3,1}$ is the subspace of prime long knots. This can be thought of
as a precise `space level' non-uniqueness result for the connect-sum
decomposition of knots, whereas Schubert's result states uniqueness on
the level of isotopy classes of knots.

There is a major conceptual gap between the Cerf--Morlet `Comparison
Theorem' and the freeness of $\K_{3,1}$ as a 2--cubes object. Getting a
better understanding of this defect was one of the primary motivations
behind this paper.

\begin{defn}\label{cubesdef}\quad
\begin{itemize}
\item A (single) little $n$--cube is
a function $L \co  \I^n \to \I^n$ such that $L=l_1\times \cdots \times l_n$
where each $l_i \co  \I \to \I$ is affine-linear and increasing ie: $l_i(t)=a_i t + b_i$ for
some $0 \leq a_i < 1$ and $b_i \in \Real$.
\item Let $\CAut_n$ denote the monoid of affine-linear
automorphisms of $\Real^n$ of the form $L=l_1 \times \cdots \times l_n$ where
$l_i \co  \Real \to \Real$ affine linear and increasing, and $L(\I^n)\subset \I^n$.
\item Given a little $n$--cube $L$ a mild abuse of notation is to
consider $L \in \CAut_n$ by taking the unique affine-linear extension of 
$L$ to $\Real^n$.
\item The space of $j$ little $n$--cubes  $\Cu_n(j)$ is the space of maps
$L \co  \sqcup_{i=1}^j \I^{n} \to \I^{n}$ such that
the restriction of $L$ to the interior of its domain is an embedding,
and the restriction of $L$ to any connected component of its domain
is a little $n$--cube. Given $L \in \Cu_{n}(j)$ let $L_i$ denote the restriction 
of $L$ to the $i$th copy of $\I^n$. By convention $\Cu_n(0)$ is
taken to be a point. This makes the union $\sqcup_{j=0}^\infty \Cu_n(j)$
into an operad, called the operad of little $n$--cubes $\Cu_n$ (see May
\cite{May1}).
\item There is an action of $\CAut_n$ on $\EK{n,M}$ given by
$$\mu \co  \CAut_n \times \EC{n}{M}  \to \EC{n}{M}$$
$$ \mu(L,f) = (L \times \Id_M)\circ f \circ (L^{-1} \times \Id_M) $$
In the above formula, $L^{-1}$ is the inverse of $L$ in the group
of affine-linear isomorphisms of $\Real^n$.
The above action is denoted $\mu(L,f)=L.f$. There is an
action of $\CAut_j$ on $\K_{n,j}$ defined essentially the same way.
\end{itemize}
\end{defn}

An action of the operad of $j$--cubes on both $\K_{n,j}$ and 
$\EK{j,M}$ where the associated multiplication on $\pi_0 \K_{n,j}$ is 
the connect-sum operation, is given next.

\begin{defn}\label{obvious}
$k_i \co  \Cu_j(i) \times \left(\K_{n,j}\right)^i \to \K_{n,j}$,
$k_i \co  \Cu_j(i) \times \EK{j,M}^i \to \EK{j,M}$
is defined by the rule
$k_i(L_1,\ldots,L_i,f_1,\ldots,f_i) = L_1.f_1 \circ \cdots \circ L_i.f_i$.
In the case of the space $\K_{n,j}$, given $f,g \in \K_{n,j}$ with
disjoint support, $f \circ g$ is defined so that 
$$f\circ g(x) = \left\{ 
\begin{array}{ll}
f(x) & \text{ if } f(x)\neq x \\
g(x) & \text{ if } \text{ otherwise. }\\
\end{array}
\right.$$
\end{defn}

\fullref{littlecdef} extends the $j$--cubes action on $\EK{j,M}$ to a 
$(j{+}1)$--cubes action.

\begin{defn}\label{littlecdef}\quad
\begin{itemize}
\item Given $j$ little $(n{+}1)$--cubes, $L=(L_1,\ldots, L_j)\in \Cu_{n+1}(j)$
define the $j$--tuple of (non-disjoint) little $n$--cubes
$L^\pi = (L_1^\pi,\ldots, L_j^\pi)$ 
by the rule $L_i^\pi =l_{i,1} \times \cdots \times l_{i,n}$ where 
$L_i= l_{i,1} \times \cdots \times l_{i,n+1}$. See \fullref{fig5}.
Similarly define
$L^t \in \I^j$ by $L^t=(L_1^t,\ldots, L_j^t)$ where $L_i^t = l_{i,n+1}(-1)$.
\begin{figure}\label{fig5}
\centering
\labellist\small
\pinlabel {$\{0\}^n\times\Real$} [t] at 190 495
\pinlabel {$L^t$} [l] at 283 271
\pinlabel {$L^{\pi}$} [tl] at 128 162
\pinlabel {$\Real^n\times\{0\}$} at 520 65
\pinlabel {$L$} [l] at 80 270
\endlabellist
\includegraphics[width=5cm]{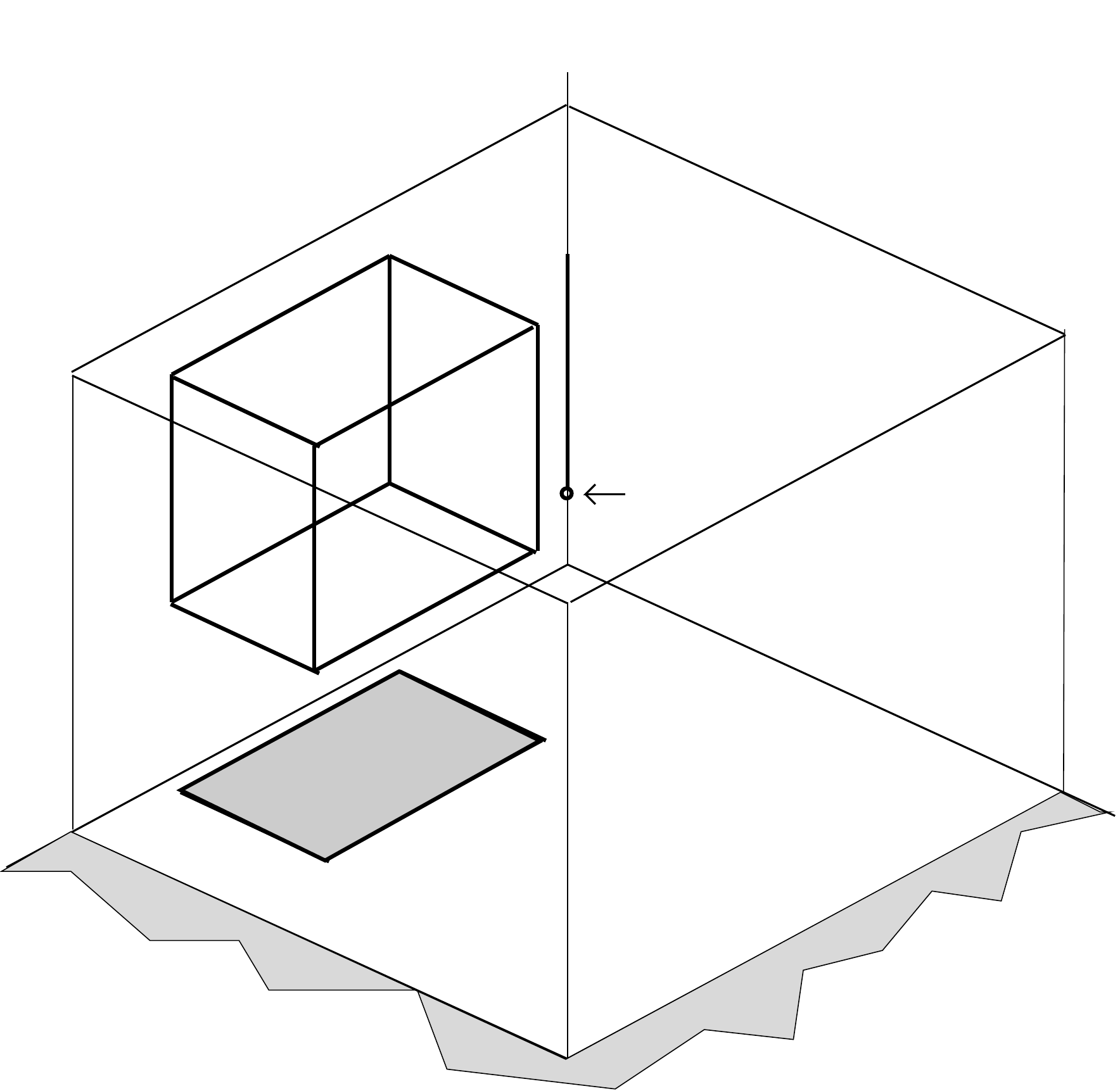}
\caption{}
\end{figure}
\item The action of the operad of little $(n{+}1)$--cubes on the space
$\EK{n,M}$ is given by
the maps $\kappa_j \co  \Cu_{n+1}(j) \times \EK{n,M}^j \to \EK{n,M}$
for $j \in \{1,2,\ldots\}$ defined by
$$\kappa_j(L_1,\ldots,L_j,f_1,\ldots,f_j) =
  L^\pi_{\sigma(1)}.f_{\sigma(1)}\circ L^\pi_{\sigma(2)}.f_{\sigma(2)}\circ\cdots \circ L^\pi_{\sigma(j)}.f_{\sigma(j)}$$
where $\sigma \co  \{1,\ldots, j\} \to \{1,\ldots,j\}$ is any permutation
such that $L^t_{\sigma(1)} \leq L^t_{\sigma(2)} \leq \cdots 
\leq L^t_{\sigma(j)}$. See \fullref{fig6}.
The map $\kappa_0 \co  \Cu_{n+1}(0) \times \EK{n,M}^0 \to \EK{n,M}$
is the inclusion of a point $*$ in $\EK{n,M}$, defined
so that $\kappa_0(*) = \Id_{\Real^n \times M}$.
\end{itemize}
\end{defn}

\begin{thm}[Budney \cite{Cube}]\label{littlecthm}
For any compact manifold $M$ and any integer $n \geq 0$
the maps $\kappa_j$ for $j \in \{0, 1,2,\ldots\}$ define an action of the
operad of little $(n{+}1)$--cubes on $\EK{n,M}$.
\end{thm}

\begin{eg}
\begin{figure}\label{fig6}
\centering
\labellist\small
\pinlabel {\tiny$L_1$} at 375 1080
\pinlabel {\tiny$L_{2}$} at 385 1375
\pinlabel {\tiny$L_3$} at 770 1240
\pinlabel {$1$} [t] at 1720 1050
\pinlabel {$-1$} [t] at 1142 1050
\pinlabel {$1$} [t] at 2520 1050
\pinlabel {$-1$} [t] at 1940 1050
\pinlabel {$1$} [t] at 3320 1050
\pinlabel {$-1$} [t] at 2742 1050
\pinlabel {$1$} [t] at 3085 40
\pinlabel {$-1$} [t] at 533 40
\pinlabel {,} at 1024 1067
\pinlabel {,} at 1827 1055
\pinlabel {,} at 2622 1065
\pinlabel {$\kappa_3$} at 1942 832
\pinlabel {$L_3^t$} [r] at 100 1135
\pinlabel {$L_2^t$} [r] at 100 1325
\pinlabel {$L_1^t$} [r] at 100 955
\pinlabel {$f_1$} at 1390 1443
\pinlabel {$f_2$} at 2183 1443
\pinlabel {$f_3$} at 2982 1443
\endlabellist
\includegraphics[width=\textwidth]{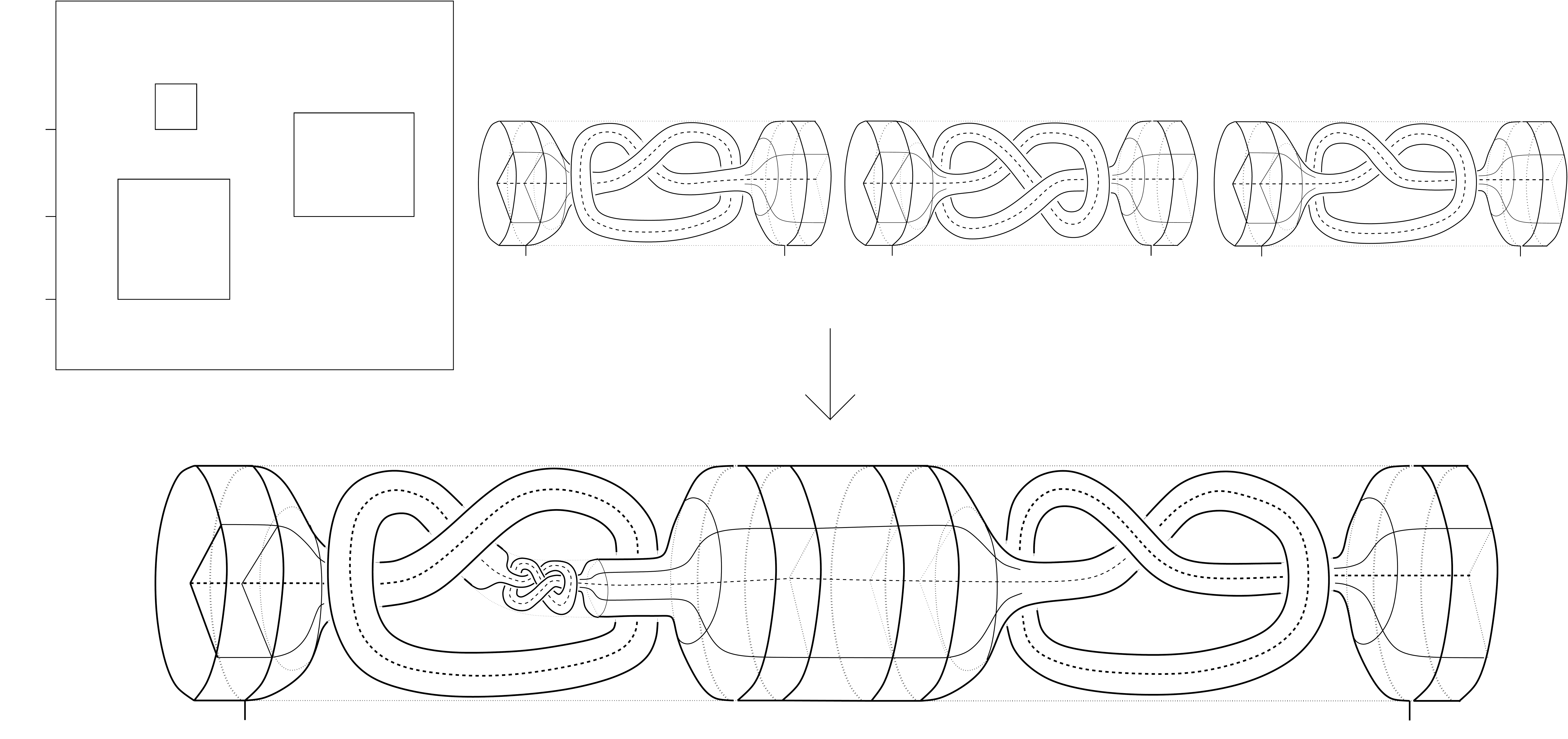}
\caption{}
\end{figure}
$L^t_1 < L^t_3 < L^t_2$ so $\sigma=(23)$ and
$\kappa_3(L_1,L_2,L_3,f_1,f_2,f_3)=L^\pi_1.f_1\circ L^\pi_3.f_3 \circ
L^\pi_2.f_2$, which explains the figure--8 knot being `inside' of the trefoil on the
left hand side of the picture. 
\end{eg}

In the definition of $\EK{n,M}$, if
one replaces the condition that the support of $f$ is contained in
$\I^n \times M$ with it being contained in $D^n \times M$ one obtains 
a homotopy-equivalent space $\ED{n,M}$. By a similar construction to
\fullref{littlecdef}, one also obtains an action 
of the operad of unframed little $(n{+}1)$--discs on $\ED{n,M}$.  
Since $\pi_0 \K_{n,j}$ is a group for $n-j>2$, $\EK{j,D^{n-j}}$ an
$(n{+}1)$--fold 
loop space. Next is a construction of analogous operad actions on the spaces 
$\PK{n,M}$. 

\begin{defn}\label{heccubes}
$\kappa_j \co  \Cu_{n}(j) \times \PK{n,M}^j \to \PK{n,M}$
for $j \in \{1,2,\ldots\}$ is defined by
$$\kappa_j(L_1,\ldots,L_j,f_1,\ldots,f_j) =
  L_{\sigma(1)}.f_{\sigma(1)}\circ L_{\sigma(2)}.f_{\sigma(2)}\circ\cdots \circ L_{\sigma(j)}.f_{\sigma(j)}$$
where $\sigma \co  \{1,\ldots, j\} \to \{1,\ldots,j\}$ is any permutation
such that $L^t_{\sigma(1)} \leq L^t_{\sigma(2)} \leq \cdots \leq L^t_{\sigma(j)}$.
\end{defn}

\begin{prop}\label{action_make}
The maps $\kappa_*$ define an action of the operad of little
$n$--cubes on $\PK{n,M}$. 
\end{prop}
\begin{proof}
There are three axioms to verify.

\begin{description}
\item[Identity] Let $\Id_{\I^n}$ be the identity $n$--cube, then
$\kappa_1(\Id_{\I^n},f)=\Id_{\I^n}.f = f$ by design.

\item[Symmetry] We need to verify that $\kappa_n(L.\alpha,f.\alpha)=
\kappa_n(L,f)$, for $\alpha \in \Sigma_n$.

Let
$$\kappa_j(L,f) =
  L_{\sigma(1)}.f_{\sigma(1)}\circ L_{\sigma(2)}.f_{\sigma(2)}\circ\cdots \circ L_{\sigma(j)}.f_{\sigma(j)}$$
and
$$\kappa_j(L.\alpha,f.\alpha) =
  L_{\alpha\sigma'(1)}.f_{\alpha\sigma'(1)}\circ L_{\alpha\sigma'(2)}.f_{\alpha\sigma'(2)}\circ\cdots \circ L_{\alpha\sigma'(j)}.f_{\alpha\sigma'(j)}$$
where $\sigma, \sigma' \in S_n$ satisfy
$L^t_{\sigma(1)} \leq \cdots \leq L^t_{\sigma(n)}$ and 
$L^t_{\alpha\sigma'(1)} \leq \cdots \leq L^t_{\alpha\sigma'(n)}$.  
Up to the ambiguity in our choice of $\sigma$ and $\sigma'$
one can assume $\sigma' = \alpha^{-1} \sigma$, giving
the result.

\item[Associativity] We need to verify the diagram below commutes:
\begin{small}
$$
\xymatrix{
\Cu_{n}(m){\times}\bigl( 
\Cu_{n}(j_1){\times}\PK{n,M}^{j_1}{\times}\cdots{\times}
\Cu_{n}(j_m){\times}
\PK{n,M}^{j_m}\bigr) \ar[r] \ar[d] &
  \Cu_{n}(m){\times}\PK{n,M}^m \ar[d] \\
\Cu_{n}(j_1{+}\cdots{+}j_m){\times}\PK{n,M}^{j_1{+}\cdots{+}j_m} 
\ar[r] & \PK{n,M}
}
$$
\end{small}
Given something in the top-left corner, consider what it maps
to in the bottom-right corner, going around both ways. 
Either way around the diagram, one gets a composite of functions of the
form $L_i.L_{i,p}.f_{i,p}$, in some order. 
The difference in the order of composition is irrelevant
as our definition only allows functions to appear in different relative
orders if they have disjoint supports.\proved
\end{description}
\end{proof}

\begin{prop}\label{Jcubemap}
Both the fibre-inclusion and projection maps in the fibration
$$\EK{n,M} \to \PK{n,M} \to \EK{n-1,M}$$
are maps of little $n$--cubes objects. The graphing map
$$\gr_1 \co  \Omega \EK{n-1,M} \to \EK{n,M}$$
is a map of $(n{+}1)$--cubes object.
\end{prop}
\begin{proof}
The map $\PK{n,M} \to \EK{n-1,M}$ is of course restriction to the 
$\{1\}\times \Real^{n-1} \times M$ `face', followed by the natural
identification with $\Real^{n-1} \times M$. 
$$\kappa_j(L_1,\ldots,L_j,f_1,\ldots,f_j) =
  L_{\sigma(1)}.f_{\sigma(1)}\circ L_{\sigma(2)}.f_{\sigma(2)}\circ\cdots \circ L_{\sigma(j)}.f_{\sigma(j)}$$
Once restricted to $\{1\}\times \Real^{n-1} \times M$ it becomes the 
composite
$$L_{\sigma(1)}^\pi.f_{\sigma(1)|\{1\}\times \Real^{n-1} \times M}\circ L_{\sigma(2)}^\pi.f_{\sigma(2)|\{1\}\times \Real^{n-1} \times M}\circ\cdots \circ L_{\sigma(j)}^\pi.f_{\sigma(j)|\{1\}\times \Real^{n-1} \times M}$$
which is precisely
$$\kappa_j(L_1,\ldots,L_j,f_{1|\{1\}\times \Real^{n-1} \times M},\ldots,
 f_{j|\{1\}\times \Real^{n-1} \times M}).$$

Consider the $(n{+}1)$--cubes action on $\Omega \EK{n-1,M}$.
Given $i$ little $(n{+}1)$--cubes $L=(L_1,\ldots,L_i)$ let
$L^\alpha = (L_1^\alpha,\ldots,L_i^\alpha) \in \Cu_1(1)^i$ be 
the projection on the 1st coordinate, and let
$L^\beta = (L_1^\beta,\ldots,L_i^\beta) \in \Cu_j(1)^i$ be
their projections on the remaining $n$ coordinates. 
The $(n{+}1)$--cubes action on 
$\Omega \EK{n-1,M}$ is given by $\kappa'$ defined below:
\begin{multline*}
\kappa'_i(L_1,\ldots,L_i,f_1,\ldots,f_i) :=  
  \kappa_i(L^\beta_1,\ldots,L^\beta_i,L^\alpha_1.f_1,\ldots,L^\alpha_i.f_i) \\
= L^{\beta \pi}_{\sigma(1)}.L^\alpha_{\sigma(1)}.f_{\sigma(1)}\circ 
  L^{\beta \pi}_{\sigma(2)}.L^\alpha_{\sigma(2)}.f_{\sigma(2)}\circ\cdots\circ 
  L^{\beta \pi}_{\sigma(i)}.L^\alpha_{\sigma(i)}.f_{\sigma(i)} 
\end{multline*}
$L^\alpha_i.f_i$ is the $\Cu_1$--action on $\Omega \EK{n-1,M}$ (reparametrisation
in the loop-space coordinate) and $L^\beta_i$ acts on this via the 
$\Cu_n$--action on $\EK{n-1,M}$. $\sigma \in \Sigma_i$ is any permutation such 
that $L^{\beta t}_{\sigma(1)} \leq L^{\beta t}_{\sigma(2)} \leq \cdots 
\leq L^{\beta t}_{\sigma(i)}$.

Consider applying the map $\gr_1$:
$$\gr_1 \co  \Omega \EK{n-1,M} \ni F \longmapsto \left( (t_0,t,v) \longmapsto
(t_0,F(t_0)(t,v)) \right) \in \EK{n,M}$$

Observe that $\gr_1 (L^{\beta \pi}_{\sigma(p)}.L^\alpha_{\sigma(p)}.f_{\sigma(p)}) = 
L^\pi_{\sigma(p)}.\gr_1 (f_{\sigma(p)})$ thus
\begin{multline*}
\gr_1(\kappa'_i(L_1,\ldots,L_i,f_1,\ldots,f_i)) = \\
  L^\pi_{\sigma(1)}.\gr_1 (f_{\sigma(1)}) \circ
  L^\pi_{\sigma(2)}.\gr_1 (f_{\sigma(2)}) \circ \cdots \circ
  L^\pi_{\sigma(i)}.\gr_1 (f_{\sigma(i)}) = \\
\kappa_i(L_1,\ldots,L_i,\gr_1(f_1),\ldots,\gr_1(f_i))
\end{multline*}
since $\gr_1$ commutes with $\circ$.
\end{proof}


\section{Survey}\label{survey}

Much of this paper has been devoted to studying the map 
$\gr_1 \co  \Omega \K_{n-1,j-1} \to \K_{n,j}$ and the pseudoisotopy 
formalism for embedding spaces. This section is more survey in 
nature, mentioning what is known on the homotopy-type of the embedding 
spaces $\K_{n,j}$ and the properties of natural maps into and out 
of these spaces, focusing largely on the issues most closely related
to iterated loop-space structures on these spaces and $\EK{j,D^{n-j}}$. 

\fullref{nullhomotopy} is a generalisation of the classical theorem that
an embedding of $S^1$ in $S^3$ unknots in $S^4$. It is based loosely on the
argument in Rolfsen's textbook \cite{Rolfsen}. The argument itself is likely
much older.

\begin{prop}\label{nullhomotopy}
The natural inclusion $\Real^n \to \Real^{n+1}$ induces an inclusion 
$i \co  \K_{n,1} \to \K_{n+1,1}$ which is null-homotopic.
\end{prop}
\begin{proof}
Two null-homotopies of $i$ will be constructed, giving a map 
$\K_{n,1} \to \Omega \K_{n+1,1}$.  

Let $j_t \co  \K_{n,1} \to \K_{n,1}$ for $t \in \I = [-1,1]$ be defined
as
$$j_t(f)(x) = \frac{f((1+t^2)x-t^3)+(t^3,0,\ldots,0)}{1+t^2}.$$
$j_0$ is the identity, yet $j_1$ consists of knots which are standard
outside of $[0,1]$, and $j_{-1}$ consists of knots which are standard
outside of $[-1,0]$. 

Let $b \co  \Real \to \Real$ be
a $C^\infty$--smooth function with the properties that:
\begin{itemize}
\item $b(x) = 0$ for all $|x| \geq 1$. 
\item $b(x)=b(-x)$ for all $x \in \Real$.
\item $b'(x) > 0$ for all $-1<x<0$.
\end{itemize}
Let $B \co  \Real \to \Real^{n+1}$ satisfy $B(x)=(x,0,\ldots,0,b(x))$.
Let $C \co  \Real \to \Real^{n+1}$ satisfy $C(x)=(x,0,\ldots,0,0)$.

Given $f \in \K_{n,1}$, consider the function
$F \co  \I  \times \Real \to \Real^{n+1}$
defined as
$$ F_t(x) = \left\{
\begin{array}{ll}
i(j_{3t}(f))(x) & \text{for } |t| \in [0,\frac{1}{3}], x \in \Real \\
(2-3|t|)\ i\left(j_{\frac{t}{|t|}}(f)\right)(x) + (3|t|-1)B(x) & \text{for
} |t| \in \bigl[\frac{1}{3},\frac{2}{3}\bigr], x \in \Real \\
(3-3|t|)B(x) + (3|t|-2)C(x) & \text{for } |t| \in \bigl[\frac{2}{3},1\bigr], x \in \Real
\end{array}
\right.
$$
$F$, restricted to either $[0,1]\times \Real$ or $[-1,0]\times \Real$ is a 
null-homotopy of $i$.
\end{proof}

It is not known whether or not $F \co  \K_{n,1} \to \Omega \K_{n+1,1}$ is null-homotopic. 
The adjoint of $F$, $\Sigma \K_{n,1} \to \K_{n+1,1}$ is the 
direct-analogue of the `Freudenthal suspension map for configuration 
spaces' (see Cohen, Cohen and Xicot\'encatl \cite{CCX}) $\Sigma C_k \Real^n \to C_k \Real^{n+1}$
which is known to induce an isomorphism on the 1st non-trivial 
homology groups of the spaces provided $n>1$. 
But in this case, first non-trivial homology group of $\Sigma \K_{n,1}$ is in 
dimension $2n-5$, while for $\K_{n+1,1}$ it is in dimension $2n-4$. 

Using the same constructions, one can construct null-homotopies of the inclusions 
$\K_{n,j} \to \K_{n+j,j}$ for all $j > 0$.

\begin{question}\quad
\begin{itemize}
\item For each $n$ and $j$, what is the smallest $i$ such that inclusion
$\K_{n,j} \to \K_{n+i,j}$ is null-homotopic?
\item Is $F \co  \Sigma \K_{n,1} \to \K_{n+1,1}$ defined in \fullref{nullhomotopy} null-homotopic?
\item If the answer to the previous question is positive, then does $F$
have two distinct null-homotopies? Is there a `Freudenthal suspension map' 
$\Sigma^2 \K_{n,1} \to \K_{n+1,1}$ inducing an isomorphism of
$H_{2n-4} \Sigma^2 \K_{n,1}$ and $H_{2n-4} \K_{n+1,1}$? 
\end{itemize}
\end{question}

There is a `fibrewise restriction' map 
$R \co  \K_{n,j} \to \Omega \K_{n,j-1}$, thinking of $\Real^j$ as 
$\Real \times \Real^{j-1}$. 
If $2n-3j-3 \geq 0$ this map is exactly $(2n{-}3j{-}3)$--connected,
as the first non-trivial homotopy groups of the two spaces are in
different dimensions.  These maps have been studied in some detail by
Morlet and Goodwillie.  The `Morlet Disjunction Lemma' (see for example
Goodwillie \cite[page~9]{GoodD}) is a theorem on the connectivity of
this map in the context of arbitrary pseudoisotopy embedding spaces.

\begin{prop}
The map $R$ is a homotopy-equivalence $R \co  \K_{n,n} \to \Omega \K_{n,n-1}$. 
\end{prop}
\begin{proof}
There are homotopy-equivalences 
$\K_{n,n} \simeq \EK{n,*}$ and $\K_{n,n-1} \simeq \EK{n-1,\I}$
given by the fibrations in \fullref{Trivprop}. 
Restriction to $\Real^{n-1} \times \I$ gives a map
$\EK{n,*} \to \EK{n-1,\I}$ which is homotopic to a fibration,
whose fibre has the homotopy-type of $\EK{n,*}^2$. The fibre-inclusion 
map $\EK{n,*}^2 \to \EK{n,*}$ is homotopic to multiplication
in the group $\EK{n,*}$ (the homotopy is constructed via the 
$(n{+}1)$--cubes action on $\EK{n,*}$). Thus, the homotopy fibre
of the map $\EK{n,*}^2 \to \EK{n,*}$ is $\EK{n,*}$. 
By \fullref{hom-fib}, this homotopy-fibre has the 
homotopy-type of $\Omega \EK{n-1,\I}$. With some additional work, 
we can see that this homotopy-equivalence is homotopic to $R$.
\end{proof}

The above argument is a mild variant of Hatcher's arguments where he gives
various equivalent statements of the Smale conjecture \cite{Hatcher2}. A
way to look at the above proposition is that studying the homotopy-type
of the spaces $\Emb(S^{n-1},S^n)$ and $\Diff(S^n)$ ultimately
reduces to studying the homotopy-types of the spaces $\K_{n,n-1}$ and
$\K_{n,n}$.  Since $\Omega \K_{n,n-1} \simeq \K_{n,n}$, the study of
the homotopy-properties of these spaces is essentially identical modulo
$\pi_0 \K_{n,n-1} \simeq \pi_0 \Emb(S^{n-1},S^n)$. The next result
compiles the major theorems on $\pi_0 \K_{n,n-1}$.

\begin{thm}\quad
\begin{itemize}
\item If $f \co  S^{n-1} \to S^n$ is a smooth embedding, then $f(S^{n-1})$
bounds a topological disc. See Mazur \cite{Maz} and Brown \cite{Brown2}.
\item The disc $D^n$ has a unique smooth structure for $n \geq 6$, 
and $D^5$ admits a unique smooth structure which restricts to the
standard smooth structure on $\partial D^5$. See Smale \cite{Smale2}.
\item (Corollary of the above two results) If $f \co  S^{n-1} \to S^n$ is a smooth 
embedding, then $f(S^{n-1})$ bounds a smooth disc provided $n \geq 5$. Thus,
the space $\Emb(S^{n-1},S^n)/\Diff(S^{n-1})$ is connected. See Kosinski \cite{Kos} for a
modern account of the results in Smale's paper \cite{Smale2}.
\item For $n \in \{2,3\}$, $\Emb(S^{n-1},S^n)$ is known to be connected. For
$n=2$ this is the Schoenflies theorem. See Siebenmann \cite{Sieb} for a historical account. For $n=3$ it is the combination of Alexander's theorem
\cite{AL}, and Smale's theorem \cite{Smale}.
\item Whether or not $\Emb(S^3,S^4)$ is connected is called the smooth
Schoenflies problem in dimension $4$.  Scharlemann \cite{Scharl} and 
Poenaru \cite{Poenaru} have some partial results on this problem.
\end{itemize}
\end{thm}

Observe that an element of $\Emb(S^{n-1},S^n)$ is isotopic to the standard
inclusion if and only if it extends to an embedding of $D^n$ in $S^n$, 
thus the kernel of the map $\pi_0 \K_{n-1,n-1} \to \pi_0 \K_{n,n-1}$
is the image of $\pi_0 \PE_{n,n} \to \pi_0 \K_{n-1,n-1}$.
The above observation that $\pi_0 \Emb(S^{n-1},S^n)/\Diff(S^{n-1})$ 
is connected for $n \geq 5$ allows the extension of the homotopy long 
exact sequence of the fibration
$\K_{n,n} \to \PE_{n,n} \to \K_{n-1,n-1}$ from \fullref{Splitprop} 
to the `classical' sequence:
$$\cdots \to \pi_1 \K_{n-1,n-1} \to \pi_0 \K_{n,n} \to \pi_0 \PE_{n,n} \to \pi_0 \K_{n-1,n-1} \to \pi_0 \K_{n,n-1} \to 0.$$
Thus, for $n \geq 5$ $\pi_0 \K_{n,n-1}$ is isomorphic to 
the groups of homotopy $n$--spheres $\theta^n$ (see Kosinski \cite{Kos}). 
$\theta^n$ is known to be finite, and many of these groups have
been computed, for example 
$\theta^5 = 0$, 
$\theta^6 = 0$,
$\theta^7 \simeq \Zed_{28}$,
$\theta^8 \simeq \Zed_2$,
$\theta^9$ is known to have 8 elements, 
$\theta^{10}$ is known to have 6 elements,
$\theta^{11} \simeq \Zed_{992}$.

\begin{thm}[Cerf \cite{Cerf2}] $\PE_{n,n}$ is connected for $n \geq 6$. 
So there is an isomorphism of groups 
$\pi_0 \Diff(D^{n-1}) \simeq \pi_0 \Emb(S^{n-1},S^n)$
and an epimorphism $\pi_1 \Diff(D^{n-1}) \to \pi_0 \Diff(D^n)$.
\end{thm}

A metric $g$ on $S^n$ is said to be round if for any points $x,y \in S^n$
there is an isometry of $g$ carrying $x$ to $y$ which can also be chosen
to send an orthonormal basis in $T_xS^n$ to any orthonormal basis in
$T_yS^n$. Let $\Bbb M^n$ denote the space of round Riemann metrics on
$S^n$. 

\begin{prop}[Hatcher \cite{Hatcher2}]
\label{met}
$\Bbb M^n$ has the same homotopy-type as $\K_{n,n} \simeq \Diff(D^n)$. 
\end{prop}
\begin{proof}
There is a fibration $\Bbb M^n \to (0,\infty)$ given by taking the
volume of the metric. The fibre of this map is a $\Diff^+(S^n)$--homogeneous 
space, with isotropy group $\SO_{n+1}$.  \fullref{Splitprop}
tells us that $\K_{n,n} \simeq \Diff(D^n)$ is also the base-space of
such a homotopy-fibre sequence $\SO_{n+1} \to \Diff^+(S^n) \to \Diff(D^n)$.
\end{proof}

Smale \cite{Smale} and Hatcher \cite{Hatcher2} have proved that $\Diff(D^n)$ is 
contractible for $n=2$ and $n=3$ respectively.  That $\Diff(D^1)$ is contractible
follows from an averaging argument, or equivalently from the `length' classification 
of connected closed  $1$--dimensional Riemann manifolds via \fullref{met}.
The space of Riemann metrics on $S^n$ is contractible since it is an affine space,
making the homotopy-type of $\Diff(D^n)$ the complete obstruction to $\Bbb M^n$ being a 
deformation-retract of the space of all Riemann metrics on $S^n$.

$\Diff(D^n)$ is an $(n{+}1)$--fold loop space (see Budney \cite{Cube},
Morlet \cite{Mor} and Burghelea and Lashof \cite{BL}) whose 
$(n{+}1)$--fold delooping is $PL(n)/O_n$ \cite{BL,Mor}. As of yet, their does
not appear to be any direct methods of studying the homotopy-type of $PL(n)$.
In particular, essentially nothing is known about the homotopy-type of
$\Diff(D^4)$. Farrell and Hsiang computed the rational homotopy of $\Diff(D^n)$ in
a range.

\begin{thm}[Farrell and Hsiang \cite{FaHs}]
Provided $0 \leq i < \min\{\frac{n-4}{3},\frac{n-7}{2}\}$ 
$$
\pi_i \Diff(D^n)\otimes \Rat \simeq
\left\{
\begin{array}{ll}
\Rat & \text{provided } 4 | (i+1) \\
0 & \text{ otherwise } \\
\end{array}
\right.
$$
\end{thm}

The bound $i < \min\{\frac{n-4}{3},\frac{n-7}{2}\}$ is known as
Igusa's stable range \cite{IG}.  Roughly this the range where $\pi_i \PE_{n,n}$ 
can be related to $K$--theory. Antonelli, Burghelea and Kahn had shown earlier
that $H_* \Diff(D^n)$ is not finitely-generated for $n \geq 7$  \cite{ABK2}.

The spaces $\K_{n+2,n}$ are in the realm of `traditional' co-dimension 
$2$ knot theory, on which there is a plethora of literature. The majority of 
the literature focuses on issues related to isotopy classification, ie: 
$\pi_0 \K_{n+2,n}$.  Some good general references are Kawauchi \cite{Kaw},
Hillman \cite{Hill,Hill2}, Ranicki \cite{Ranicki} and Kervaire--Weber \cite{KerWeb}. 

The homotopy-type of $\K_{3,1}$ has been described,
component-by-component, as an iterated fibre bundle in the author's
article \cite{KnotSpace}, which builds on the previous works of Hatcher
\cite{Hatcher3,Hatcher4}, and the author \cite{comp_tree,Cube}.

\begin{thm}[Budney \cite{KnotSpace}]\label{k31ans}
Given a long knot $f \in \K_{3,1}$, let $\K_{3,1}(f)$ denote the path component
in $\K_{3,1}$ containing $f$.  Then $\K_{3,1}(f)$ has the homotopy-type of:
\begin{enumerate}
\item $\{*\}$ if $f$ is the unknot.
\item $S^1 \times \K_{3,1}(g)$ if $f$ is a cable of $g$. 
\item $C_n(\Real^2) \times_{\Sigma_f} \prod_{i=1}^n \K_{3,1}(f_i)$ if 
$f= f_1 \# \cdots \# f_n$ is the prime decomposition of $f$, with $n \geq 2$.
$\Sigma_f$ is the subgroup of $\Sigma_n$ corresponding to the partition
of $\{1,2,\ldots,n\}$ defined by the equivalence relation $i \sim j$ if and only
if $\K_{3,1}(f_i)=\K_{3,1}(f_j)$.
\item $S^1 \times \left( SO_2 \times_{A_f} \prod_{i=1}^n \K_{3,1}(f_i) \right)$
if $f = (f_1,\ldots,f_n)\splice L$ is hyperbolically-spliced. Here $L$ is
some hyperbolic link $L=(L_0,L_1,\ldots,L_n)$ in $S^3$ with the $L_0$ component
`long'. Define $B_L$ to be the group of orientation-preserving hyperbolic isometries
of $S^3 \setminus L$ which extend to $L$, preserving $L_0$ and its orientation. 
$B_L \to \Diff(S^3,L_0)$ is a faithful representation, giving an embedding of
$B_L$ in $\Diff(L_0)$ (thus conjugate to a subgroup of $SO_2$).  Similarly, 
there is a homomorphism $B_L \to \pi_0 \Diff(L_1 \cup \ldots \cup L_n) \equiv \Sigma_n^+$ 
the signed symmetric group of $\{1,2,\ldots,n\}$. $\Sigma_n^+$ acts on 
$\K_{3,1}^n$ by permutation of the factors and knot inversion.
Let $A_f$ be the subgroup of $B_L \subset \Sigma_n^+$ that preserves 
$\prod_{i=1}^n \K_{3,1}(f_i)$.
\end{enumerate}
\end{thm}

Case (2) above is considered to apply to torus knots -- think of a torus knot
as a cable of the unknot, thus the component of a torus knot has the homotopy-type
of $S^1$. A hyperbolic knot is thought of as a hyperbolically-spliced
knot where $L$ is a $1$--component hyperbolic link, thus such a component has
the homotopy-type of $S^1\times S^1$. Since every knot can be 
obtained from the unknot by iterated cabling, connect-sum and hyperbolic splicing 
operations \cite{comp_tree}, the above result describes the homotopy-type of
$\K_{3,1}(f)$ for any $f \in \K_{3,1}$.  To be clear, if the knot $f$ has $j$ tori in the JSJ-decomposition of its complement, to obtain an answer for the homotopy-type of 
$\K_{3,1}(f)$, one would have to apply \fullref{k31ans} $j+1$ times. 
A detailed justification for the above theorem is given in the reference 
\cite{KnotSpace}. The homotopy-equivalence in part (3) of \fullref{k31ans} is induced by
the action of the operad of $2$--cubes on $\K_{3,1}$.  Another way to state (3)
is that $\K_{3,1}$ is a free $2$--cubes object, with generating space
$\Prime \sqcup \{*\}$, $\Prime \subset \K_{3,1}$ the space of prime long knots. 
By the work of May \cite{May1}, the group-completion $\Omega B\K_{3,1}$ of
the knot space has a particularly simple structure, 
$\Omega B\K_{3,1} \simeq \Omega^2 \Sigma^2 \left( \Prime \sqcup \{*\} \right)$.
Fred Cohen and the author have used these results to compute the homology of
many components of $\K_{3,1}$ \cite{BudCoh}. In the process it became clear
that the homotopy-type and homology of $\K_{3,1}$ would likely have a more elegant 
description if one could prove that $\K_{3,1}$ had an action of the operad of 
framed little $2$--discs. 

\begin{question}
Can one define an action of the operad of framed $(n{+}1)$--discs on the
spaces $\ED{n,D^k}$, in a `natural geometric manner' similar to \fullref{littlecdef}? $\ED{n,D^k}$ refers to the comments preceding
\fullref{heccubes}.
\end{question}


The topic of $\pi_0 \K_{4,2}$ has a few new references. Carter and Saito have 
constructed an analogue of Reidermeister theory \cite{CartSeit}. Kamada has 
constructed an analogue of the Alexander--Markov theorem from dimension $3$ 
\cite{Kamada}. 

\begin{thm}[Zeeman \cite{ZeemTwist} and Litherland \cite{Lith}]
Let $g \in \Omega\K_{n+2,n}(f)$ be such that $\tilde g \in \pi_0 \Diff(\I^{n+2},f)$ 
preserves a Seifert surface for $f$. Let $G \in \pi_0 \Diff(\I^{n+2},f)$ 
denote the Gramain element (a meridional Dehn twist). If 
$k \in \Zed \setminus \{0\}$ then the complement of $\gr_1(G^kg) \in \K_{n+3,n+1}$ 
fibres over $S^1$.
\end{thm}

For $n=1$ Litherland went on to identify the fibre in several cases.
From a practical point of view, the Zeeman--Litherland theorem is a useful
tool for constructing embeddings of 3--manifolds into $S^4$, as fibres of
fibred knot complements (see Ruberman \cite{Ruberman}).  It is possible
that there are other types of Alexander--Markov theorems in dimension four.
Recently it was shown by Mozgova and the author that Litherland spinning
does not suffice \cite{budmoz}, because the Alexander polynomial provides
an obstruction to elements of $\pi_0 \K_{4,2}$ being deform-spun.

Up to a homotopy-equivalence, the spaces $\ED{j,D^{n-j}}$ and
$\EK{j,D^{n-j}}$ admit an action of the operad of framed little
$(j{+}1)$--discs, provided $n-j>2$.  This is because they are
$(j{+}1)$--fold loop spaces.  This argument does not apply when $n-j =
2$ since $\pi_0 \EK{n,D^2}$ is never a group.  This will be explained
in the next proposition.

\begin{prop}\quad
\begin{itemize}
\item $\pi_0 \K_{n+2,n}$ is not a group for all $n \geq 1$.
\item The map $\pi_0 \K_{n+1,n} \to \pi_0 \K_{n+2,n}$ induced by
inclusion $\Real^{n+1} \to \Real^{n+2}$ is injective
and maps onto the maximal subgroup of $\pi_0 \K_{n+2,n}$
provided $n \geq 4$.
\end{itemize}
\end{prop}
\begin{proof}
To prove the first point, non-invertible elements are constructed.
Start with $f_1 \in \K_{3,1}$ a trefoil knot.  Then $\pi_1 C_f$ is the
braid group on $3$ strands.  Let $g_1 = 0 \in \pi_1 \K_{3,1}(f_1)$
be the constant loop, and observe that the complement of $f_2 =
\gr_1(g_1) \in \K_{4,2}$ also has the braid group on $3$ strands as
its fundamental group. Continuing, this constructs for all $n \geq 1$
a knot $f_n \in \K_{n+2,n}$ whose complement has the braid group on $3$
strands as its fundamental group. $f_n$ is non-invertible in the monoid
$\pi_0 \K_{n+2,n}$ by Zieschang, Vogt and Coldewey \cite[Proposition
2.3.4]{ZVC}.  This is because if $h \in \K_{n+2,n}$ then the complement
of the connect-sum $f_n \# h$, $C_{f_n \# h}$ has the homotopy-type of
the union of $C_{f_n}$ and $C_h$ where $C_{f_n}$ and $C_h$ intersect
along a meridional circle, so by the canonical form for amalgamated free
products, $\pi_1 C_{f_n \# h}$ contains $\pi_1 C_{f_n}$.

By the above argument, if $f \in \pi_0 \K_{n+2,n}$ is invertible,
$\pi_1 C_f \simeq \Zed$. By a Mayer--Vietoris sequence argument, 
$H_i C_f = 0$ for all $i > 1$. Thus, $C_f$ has the homotopy-type of
a circle. By Levigne's unknotting theorem \cite{Lev} (provided $n \geq 4$)
or Wall's unknotting theorem \cite{Wall} (for $n=3$), $f$ is in the
image of $\pi_0 \K_{n+1,n}$. 

The last item to prove is that the map $\pi_0 \K_{n+1,n} \to \pi_0 \K_{n+2,n}$
is injective.  Consider $S^n \subset S^{n+1} \subset S^{n+2}$.  Let
$f \co  S^n \to S^{n+2}$ be an embedding with $f(S^n)=S^n$.  By \fullref{Splitprop} we could equivalently prove that if $f$ extends to
an embedding $F \co  D^{n+1} \to S^{n+2}$, then there is another extension of
$f$, $F' \co  D^{n+1} \to S^{n+1}$.  Identify the complement of an open 
tubular neighbourhood of $S^n$ in $S^{n+2}$ with $S^1 \times D^{n+1}$. 
Thus, $F$, if it exists, is an embedding $F \co  D^{n+1} \to S^1 \times D^{n+1}$
such that $F(\partial D^{n+1})=\{1\}\times \partial D^{n+1}$. 
By Farrell's proof of the relative Browder--Livesay--Leving--Farrell fibration 
theorem \cite{Farr}, there is a diffeomorphism 
$G \co  S^1 \times D^{n+1} \to S^1 \times D^{n+1}$ such that 
$G(F(D^{n+1})) = \{1\}\times D^{n+1}$ and $G_{|S^1 \times \partial D^{n+1}}$
is the identity on $S^1 \times \partial D^{n+1}$.  Farrell's theorem requires
$n \geq 4$. The basic idea of the proof is much like the proof of
\fullref{Embdk}, but in this case one lifts $F$ to an embedding
of $D^{n+1}$ in $\Real \times D^{n+1}$ and applies the relative H-cobordism
theorem to acquire the neccessary diffeomorphism.
\end{proof}

I would like to thank Larry Siebenmann for suggesting the 
Browder--Livesay--Leving--Farrell fibration theorem.

The above proposition implies that $\EK{n,D^2}$ is not a free
$(n{+}1)$--cubes object provided there exists exotic $(n{+}1)$--spheres, 
so no direct analogue of \cite{Cube} is true in high dimensions. Of
course, $\EK{1,D^2}$ is not a free object, either, as it splits as 
a product of $\Zed$ with the free object $\K_{3,1}$. 
One might hope that for $n>1$, $\EK{n,D^2} \simeq \K_{n+2,n}$ is closely 
related to a free $(n{+}1)$--cubes object, but there are yet further
obstructions. Provided $n \geq 3$, 
$\pi_0 \K_{n+2,n} / \pi_0 \K_{n+1,n}$ (this is the isotopy classes of
the images of the elements of $\K_{n+2,n}$) is not a free commutative monoid.
Kearton proved this in the $n=3$ case, which has since been generalised to
all $n \geq 3$. Bayer--Fluckiger went on to prove the non-existence of a `cancellation
law' ie: one can satisfy $a+b=a+c$ with $b \neq c$. See Kearton's survey \cite{Kearton}
for details.

\begin{question}\label{freeq}\quad
\begin{itemize}
\item What is the group-completion of the monoid $\pi_0 \K_{n+2,n}$? 
\item Can one characterise the monoid structure on $\pi_0 \K_{n+2,n}$
for $n \geq 2$? 
\item If $f \in \K_{n+2,n}$ is a connect-sum of two non-trivial knots,
the action of the operad of $(n{+}1)$--cubes on $\K_{n+2,n}$ gives a map 
$S^n \to \K_{n+2,n}(f)$. Is this map a non-trivial element of
$\pi_n \K_{n+2,n}(f)$? 
\end{itemize}
\end{question}

For the last of the above questions, a theorem of Swarup's \cite{Swarup} is
relevant.  He proves that if $C_f$ is the complement of a non-trivial
co-dimension two knot $f \in \K_{n+2,n}$ with $n > 2$ then the knot longitude
is a non-trivial element of $\pi_n C_f$. 

The remainder of the survey will focus on the high co-dimension
case: $\K_{n,j}$ for $n-j>2$. For references, Adachi's survey has been 
around for a few years \cite{AD}. It focuses on topics such as the Whitney trick,
and the Smale--Hirsch immersion theorem. Skopenkov has a recent survey article 
\cite{Skop} which is concerned with $\pi_0 \K_{n,j}$.   
Goodwillie, Klein and Weiss have recently put together a survey of 
what is known about embedding spaces from the point of view of disjunction 
\cite{GoodDis}. 

There have been computations of some of the groups $\pi_0 \K_{n,j}$. 
From \fullref{Connectivityknj}, the first non-trivial homotopy-group 
of $\K_{n,j}$ is in dimension $2n-3j-3$ (provided $2n-3j-3 \geq 0$).  Along the 
$2n-3j-3=0$ line there is $\pi_0 \K_{3,1}$ which is the free commutative monoid on 
$\pi_0 \Prime$, the isotopy-classes of prime long knots (see Schubert
\cite{Sch}). 
Provided $j>1$ and $2n-3j-3=0$, there are Haefliger's computations \cite{Haefliger2}:
$$
\pi_0 \K_{n,j} \simeq
\left\{
\begin{array}{ll}
\Zed & j \equiv 3(\text{mod }4) \\
\Zed_2 & j \equiv 1(\text{mod }4) \\
\end{array}
\right.
$$ 
The generator being Haefliger's Borromean rings construction \cite{Haefliger1},
also sometimes called the `trefoil' \cite{Skop}. 
The generator has also been described (\fullref{resgen}) as an iterated graphing
construction applied to $r$, the resolution of an immersion of $\Real$ in 
Euclidean space, corresponding to the $\bigotimes$ chord-diagram (see
Cattaneo, Cotta-Ramusino and Longini \cite{Cat1}).
More recently, another spinning construction involving
$r$ has recently been developed by Roseman and Takase \cite{Roseman}.

The work of Haefliger \cite{Haefliger2}, Milgram \cite{Mil}, 
Kreck and Skopenkov \cite{KrSk} gives $\pi_0 \K_{n,j}$ along
the $n-j>2$ part of the $2n-3j-3=-1$ line. Their computations
are:
$$
\pi_0 \K_{n,j} \simeq
\left\{
\begin{array}{ll}
0 &  j \equiv 2 \text{ or } 6 (\text{mod } 4) \\
\Zed_{12} & (n,j) = (7,4) \\
\Zed_4 & j \equiv  4 (\text{mod } 8), j \geq 12 \\
\Zed_2\oplus \Zed_2 & j \equiv 0 (\text{mod } 8)
\end{array}
\right.
$$ 
The above results give the next corollary as a direct analogue to 
\fullref{Connectivityknj}.
\begin{cor}\quad
\begin{itemize}
\item $\pi_{6n}\K_{3n+4,2}$ is non-trivial and has $\Zed_2\oplus \Zed_2$ as a quotient
for all $n \geq 1$.
\item $\pi_{6n+2}\K_{3n+5,2}$ is non-trivial and has $\Zed_4$ as a quotient for all $n \geq 0$
($\Zed_{12}$ for $n=0$).
\end{itemize}
\end{cor}

\begin{question}
What is the structure of the groups $\pi_2 \K_{5,2}$ and $\pi_6 \K_{7,2}$.
Further, find explicit geometric representatives for the embeddings, in analogy to \fullref{resgen}.
\end{question}

The technique of Haefliger \cite{Haefliger2} involves
two main steps. The first step is the construction of
an isomorphism $\pi_0 \K_{n,j} \simeq C^{n-j}_j$
where $C^{n-j}_j$ is the group of concordance classes of embeddings of
$S^j$ in $S^{n}$. This step is formally analogous to \fullref{Embdk}.
Using a Thom-type construction, Haefliger constructs
an isomorphism between $C^{n-j}_j$ and a multi-relative homotopy group
$C^n_j \simeq \pi_{j+1}(G;SO,G_{n-j})$ where
$SO = \smash{\varinjlim}\left(\SO_1 \to \SO_2 \to \SO_3 \to \cdots \right)$ 
is the stable special-orthogonal group, $G_n$ is the space of
degree $1$ self-maps of $S^{n-1}$, with $G$ the analogous
stable object, defined via suspensions 
$G=\smash{\varinjlim}\left(G_1 \to G_2 \to G_3 \to \cdots \right)$.
This reduces the computation of $\pi_0 \K_{n,j}$ to rather
traditional difficult problems common to surgery theory \cite{Ranicki}: 
homotopy groups of spheres and orthogonal groups.

Takase \cite{Tak} has recently proved that
any embedding of $S^{4k-1} \to S^{6k}$ can be extended to
an embedding of $(S^{2k}\times S^{2k}) \setminus D^{4k} \to S^{6k}$. 
Takase gives a rather explicit formula for determining
the isotopy class of an element of $\Emb(S^{4k-1},S^{6k})$
that simplifies Haefliger's characteristic class computations
\cite{Haefliger1}.

The work of Volic, Lambrechts and Turchin \cite{VLT} gives the
homology $H_* (\K_{n,1};\Rational)$ for $n \geq 4$ as the homology
of a differential graded algebra, by showing the collapse of the
rational Vassiliev spectral sequence. Turchin has found a Poisson
algebra structure for this DGA \cite{Tur3,Tur2}, which motivated
the author's construction of the $2$--cubes action on $\K_{3,1}$.
Salvatore \cite{Salvatore}, building on the work of Sinha \cite{Dev2} has
recently constructed a $2$--cubes action on $\K_{n,1}$ for $n \geq 4$.
The structure of $\K_{n,1}$ and $\EK{1,D^{n-1}}$ as $2$--cubes objects
for $n \geq 4$ remains mysterious.  One would hope that constructions
having the flavour of Mostovoy's \cite{Mostovoy} `short rope' spaces,
or Anderson and Hsiang's `bounded embedding spaces' \cite{AH} could give
useful geometric models that one could use to get homotopy-theoretic
information on $B^j \K_{n,j}$, $B^2 \K_{n,1}$, $B^{j+1} \EK{j,M}$. Not
only is there a lack of proofs that these spaces are the appropriate
iterated classifying spaces, but, even if they were, its not clear how
one could use such results to study the spaces $\K_{n,j}$.

\bibliographystyle{gtart}
\bibliography{link}

\end{document}